\documentclass[12pt,a4paper]{article}
\usepackage{amsmath,amsfonts,amssymb,amscd}

\usepackage{graphicx}
\usepackage{amsfonts,amssymb,amscd,amsmath}

\usepackage{indentfirst}

\usepackage{tikz}
\usetikzlibrary{matrix}

\def\H{\operatorname{H}}

\def\Ker{\operatorname{Ker}}
\def\Ad{\operatorname{Ad}}
\def\ad{\operatorname{ad}}

\def\GL{\operatorname{GL}}

\def\Im{\operatorname{Im}}
\def\pt{\operatorname{pt}}

\newcounter{th}
\def\t{\refstepcounter{th}{\bf \noindent{Theorem} \arabic{th}. }}

\newcounter{prop}
\def\prop{\refstepcounter{prop}{\bf \noindent{Proposition} \arabic{prop}. }}

\newcounter{lem}
\def\lem{\refstepcounter{lem}{\bf \noindent{Lemma} \arabic{lem}. }}

\newcounter{de}
\def\de{\refstepcounter{de}{\bf \noindent{Definition} \arabic{de}. }}

\newcounter{ex}

\begin{document}

\begin{center}

{\LARGE{\bf Vector fields on $\mathfrak{gl}_{m|n}(\mathbb C)$-flag supermanifolds}}\footnote[1]{Supported by Max Planck Institute for Mathematics, Bonn.}
\bigskip

{\bf E.G.Vishnyakova}\\[0.3cm]
\end{center}

\bigskip

\begin{abstract}
	The main result of this paper is the computation of the Lie superalgebras of holomorphic
vector fields on  complex flag supermanifolds, introduced by Yu.I.~Manin. We prove  that with several exceptions any holomorphic vector field is fundamental with respect to the
natural action of the Lie superalgebra $\mathfrak {gl}_{m|n}(\mathbb C)$.

\end{abstract}

\bigskip

\section{Introduction} 
It is a classical result that all holomorphic vector fields on a flag manifold in $\mathbb C^n$ are fundamental for the natural action of the general linear Lie group $\GL_{n}(\mathbb C)$. More precisely the Lie algebra of holomorphic vector fields on a flag manifold is isomorphic to $\mathfrak{pgl}_{n}(\mathbb C)$. Similar statement holds with some exceptions for flag manifolds that are isotropic with respect to a non-degenerate symmetric or skew-symmetric bilinear form in $\mathbb C^n$. These results were obtained by A.L.~Onishchik in 1959, see example \cite{ADima} for details.

In \cite{Man} Yu.I.~Manin constructed four series of complex compact homogeneous supermanifolds corresponding to four series of classical linear Lie superalgebras: $\mathfrak{gl}_{m|n}(\mathbb C)$, $\mathfrak{osp}_{m|n}(\mathbb C)$, $\pi\mathfrak{sp}_{n}(\mathbb C)$ and $\mathfrak{q}_{n}(\mathbb C)$, see \cite{Kac} for precise definitions. The present paper is devoted to the calculation of the Lie superalgebras of holomorphic vector fields on complex flag supermanifolds corresponding to the Lie superalgebra $\mathfrak{gl}_{m|n}(\mathbb C)$. 
 It turns out that under some restrictions on the flag type all global holomorphic vector fields are fundamental with respect to the natural action of the Lie superalgebra $\mathfrak {gl}_{m|n}(\mathbb C)$. In case of super-Grassmannians the similar result  was obtained in \cite{onigl}.

 In the present paper we study flag supermanifold $\mathbf
 F_{k|l}^{m|n}$ of type $k|l$ in the vector superspace $\mathbb C^{m|n}$. Here we put $k=(k_1,\ldots,k_r)$ and
 $l=(l_1,\ldots,l_r)$ such that
 \begin{equation}\label{eq conditions on k_i and l_j}
 \begin{split}
 0\le k_r\le\ldots\le k_1\le m&,\quad 0\le l_r\ldots\le l_1\le n\quad \text{and}\\
 \quad 0 <
 k_r +l_r <&\ldots < k_1 + l_1 < m+n.
 \end{split}
 \end{equation} 
 The number $r$ is called the {\it length} of $\mathbf
 F_{k|l}^{m|n}$.  The idea of the proof is to use results of  \cite{onigl} and the following fact. 
 For $r > 1$ the supermanifold $\mathbf
 F_{k|l}^{m|n}$ is the total space of a holomorphic superbundle with base space isomorphic to the super-Grassmannian $\mathbf
 F_{k_1|l_1}^{m|n}$ and the fiber isomorphic to a flag supermanifold of length $r-1$. The projection of this superbundle is equivariant with respect to the natural actions of the Lie supergroup $\GL_{m|n}(\mathbb C)$ on the total space and base space of $\mathbf F_{k|l}^{m|n}$.

 Let $p :\mathcal M \to \mathcal B$ be a
 morphism of supermanifolds. A vector field $v$ defined on $\mathcal M$ is said to be {\it projectable} with
 respect to $p$ if there is a vector field $v_1$  on  $\mathcal B$ such that 
 $$
 p^*(v_1(f)) = v(p^*(f))
 $$ 
 for any $f\in \mathcal O_{\mathcal B}$. A vector field $v$ on $\mathcal M$ is
called {\it vertical} if it is projected to $0$. If $p$ is a projection of a superbundle, then
 every projectable vector field $v$ is projected to a unique vector field $v_1$.  In \cite{Bash} the following statement was proven. If $p :\mathcal M \to \mathcal B$ is the projection of a superbundle with fibre $\mathcal S$ with  $\mathcal O_{\mathcal S}(\mathcal S_0)= \mathbb C$, this is any global holomorphic function on $\mathcal S$ is constant, then every vector field on $\mathcal M$ is projectable with respect to $p$.
 Denote by $\mathfrak v (\mathcal M)$ the Lie superalgebra of holomorphic vector fields on $\mathcal M$. If $\mathcal O_{\mathcal S}(\mathcal S_0)= \mathbb C$, we have a map  
 $$
 \mathcal P:\mathfrak v (\mathcal M) \to \mathfrak v (\mathcal B).
 $$
 This map is a Lie superalgebra homomorphism, and its kernel $\Ker \mathcal P$ is the set of all vertical vector fields.

Consider the superbundle $\mathbf F_{k|l}^{m|n}$. The space of global holomorphic functions $\mathcal O_{\mathcal S}(\mathcal S_0)$ was computed in \cite{Viholom}. It was shown that $\mathcal O_{\mathcal S}(\mathcal S_0)= \mathbb C$ under some restrictions on the flag type $k|l$. Therefore, in general all holomorphic vector fields on $\mathcal M$ are projectable to the super-Grassmannian $\mathcal B= \mathbf F_{k_1|l_1}^{m|n}$ and we have the following homomorphism of Lie superalgebras
 $$
 \mathcal P:\mathfrak v (\mathbf F_{k|l}^{m|n}) \to \mathfrak v (\mathbf F_{k_1|l_1}^{m|n}).
 $$
From the equivariance of $p$ with
respect to the actions of $\GL_{m|n}(\mathbb C)$ it follows that the natural Lie algebra homomorphisms
$$
\mu: \mathfrak{gl}_{m|n}(\mathbb C)\to  \mathfrak v (\mathbf F_{k|l}^{m|n})\quad \text{and} \quad
\mu_{\mathcal B}: \mathfrak{gl}_{m|n}(\mathbb C)\to  \mathfrak v (\mathbf F_{k_1|l_1}^{m|n})
$$ 
satisfy the relation $\mu_{\mathcal B}= \mathcal P\circ \mu$. Assuming that  the homomorphism $\mu_{\mathcal B}$ is
surjective, in other words assuming that 
$$
\mathfrak v (\mathbf F_{k_1|l_1}^{m|n})\simeq \mathfrak{pgl}_{m|n}(\mathbb C),
$$ 
we see that $\mathcal P$ is also surjective. The main goal of this paper is prove that $\mathcal P$ is injective. Then $\mathcal P$ is invertible and we have 
$$
\mu = \mathcal P^{-1}\circ\mu_{\mathcal B}.
$$
Therefore, 
$$
\mathfrak v (\mathbf F_{k|l}^{m|n})\simeq \mathfrak{pgl}_{m|n}(\mathbb C).
$$

The main result of this paper was announced in \cite{Vivector} in case $0< k_r<\ldots< k_1< m$ and $0< l_r\ldots< l_1< n$ and the idea of the proof was sketched in \cite{Vi} also in this case. Here we give the proof in general case. Our main result is the following.

\medskip

\noindent{\bf Theorem.} {\sl Assume that $r>1$ and that we have the following restrictions on the flag type: 
	\begin{align*}
	&(k_i,l_i)\ne (k_{i-1},0),\,(0,l_{i-1}),\,\, i\geq 2;\\
	&(k_{i-1},k_i|l_{i-1},l_i)\ne (1,0|l_{i-1},l_{i-1}-1),\, 	(1,1|l_{i-1},1),\,\, i\geq 1;\\
	&(k_{i-1},k_i|l_{i-1},l_i)\ne (k_{i-1},k_{i-1}-1|1,0),\, 	(k_{i-1},1|1,1),\,\, i\geq 1;\\
	&k|l\ne 	(0,\ldots,0| n,l_2,\ldots,l_r),\,\, k|l\ne 	( m,k_2,\ldots,k_r|0,\ldots,0).
	\end{align*}
	Then
	$$
	\mathfrak{v}(\mathbf{F}^{m|n}_{k|l})\simeq \mathfrak {pgl}_{m|n}(\mathbb
	C).
	$$
	If $k|l= (0,\ldots,0| n,l_2,\ldots,l_r)$ or $k|l= 	( m,k_2,\ldots,k_r|0,\ldots,0)$, then
	$$
	\mathfrak{v}(\mathbf{F}^{m|n}_{k|l})\simeq W_{mn}\subset\!\!\!\!\!\!+ (\bigwedge(\xi_1,\ldots,\xi_{mn})\otimes\mathfrak {pgl}_{n}(\mathbb
	C)),
	$$
	where $W_{mn}=\operatorname{Der}\bigwedge(\xi_1,\ldots,\xi_{mn})$.
	
}

\section{Preliminaries}

\subsection{Flag supermanifolds}

We will use the word ``supermanifold'' in the sense of Berezin and Leites \cite{BL}. Throughout we will restrict our attention to the complex-analytic version of the theory. Recall that a {\it complex-analytic superdomain of dimension $s|t$} is a $\mathbb{Z}_2$-graded locally
ringed space of the form 
$$
\mathcal U = (\mathcal U_0, \mathcal{F}_{\mathcal U_0} \otimes_{\mathbb C} \bigwedge(t) ),
$$
where $\mathcal{F}_{\mathcal U_0}$ is the sheaf of holomorphic functions on an open set $\mathcal U_0\subset \mathbb{C}^s$ and
$ \bigwedge(t)$ is the Grassmann algebra with $t$ generators.
 A {\it complex-analytic supermanifold} of dimension $s|t$ is a $\mathbb{Z}_2$-graded locally ringed space that
is locally isomorphic to a complex-analytic superdomain of dimension $s|t$. We will denote a supermanifold by $\mathcal{M} = (\mathcal{M}_0,{\mathcal O}_{\mathcal{M}})$, where $\mathcal{M}_0$ is the underlying complex-analytic manifold and ${\mathcal O}_{\mathcal{M}}$ is the structure sheaf.

Let us give an explicite description of a flag supermanifold in terms of charts and local coordinates (see also \cite{Man,Viholom}). Let us take two non-negative integers $m,n\in\mathbb Z$ and two sets of non-negative integers 
$$
k=(k_1,\ldots,k_r)\quad  \text{and} \quad  l=(l_1,\ldots,l_r)
$$
such that (\ref{eq conditions on k_i and l_j}) holds. 
The underlying space of the supermanifold
$\mathbf F_{k|l}^{m|n}$ is the product  $\mathbf
F^{m}_{k}\times \mathbf F^{n}_{l}$ of two flag manifolds of types $k=(k_1,\ldots,k_r)$ and $l=(l_1,\ldots,l_r)$  in the vector spaces $\mathbb C^m$ and $\mathbb C^n$, respectively. Let us fix two subsets 
$$
I_{s\bar
0}\subset\{1,\ldots,k_{s-1}\}\quad \text{and} \quad I_{s\bar 1}\subset\{1,\ldots,l_{s-1}\},
$$
 where $k_0=m$ and $l_0=n$, such that
$|I_{s\bar 0}| = k_s,$ and $|I_{s\bar 1}| = l_s$,  for any $s = 1,\ldots,r$. We put
$I_s=(I_{s\bar 0},I_{s\bar 1})$ and $I = (I_1,\ldots,I_r)$. 
We assign the following $(k_{s-1} + l_{s-1})\times (k_s + l_s)$-matrix 
\begin{equation}\label{eq local chart general}
Z_{I_s}
=
\left(
\begin{array}{cc}
X_s & \Xi_s\\
\H_s & Y_s \end{array} \right), \ \ s=1,\dots,r,
\end{equation}
to any $I_s$. Here we assume that 
$$
X_s=(x^s_{ij})\in \operatorname{Mat}_{k_{s-1}\times k_s}(\mathbb C),\quad
Y_s=(y^s_{ij})\in \operatorname{Mat}_{l_{s-1}\times l_s}(\mathbb C),
$$
are even elements and elements of the matrices 
$\Xi_s=(\xi^s_{ij})$, $\H_s=(\eta^s_{ij})$ are odd. We also assume that $Z_{I_s}$ contains the identity submatrix $E_{k_s+l_s}$ of size $(k_s+l_s)\times (k_s+l_s)$ in the lines with numbers $i\in I_{s\bar 0}$ and $k_{s-1} + i,\; i\in I_{s\bar 1}$. For example in case 
\begin{align*}
I_{s\bar 0}=\{k_{s-1}-k_s+1,\ldots, k_{s-1}\},\quad
 I_{s\bar 1}=\{l_{s-1}-l_s+1,\ldots, l_{s-1}\}
\end{align*}
the matrix $Z_{I_s}$ has the following form:
$$
Z_{I_1} =\left(
\begin{array}{cc}
X_s&\Xi_s\\
E_{k_s}&0\\
\H_s&Y_s\\0&E_{l_s}\end{array} \right).
$$
Here $E_{q}$ is the identity matrix of size $q\times q$. For simplisity of notation we use here the same letters $X_s$, $Y_s$, $\Xi_s$ and $\H_s$ as in (\ref{eq local chart general}).

We see that the sets $I_{\bar 0} =(I_{1\bar 0},\ldots,I_{r\bar 0})$ and $I_{\bar
1}= (I_{1\bar 1},\ldots,I_{r\bar 1})$ determine the charts $U_{I_{\bar 0}}$ and $V_{I_{\bar 1}}$ on the flag manifolds $\mathbf F^{m}_{k}$ and $\mathbf F^{n}_{l}$, respectively. We can take the non-trivial elements (i.e., those are not contained in the identity submatrix) from
 $X_s$ and $Y_s$ as local coordinates in $U_{I_{\bar 0}}$ and $U_{I_{\bar 1}}$, respectively. Summing up, we have defined the following atlas on $\mathbf F^{m}_{k} \times\mathbf
 F^{n}_{l}$:
$$
\{U_I = U_{I_{\bar 0}}\times
U_{I_{\bar 1}}\} 
$$
 with charts are parametrized by $I = (I_s)$. The sets $I_{\bar 0}$ and $I_{\bar 1}$ also determine the superdomain $\mathcal U_I$ with underlying space $U_I$ and with even and odd coordinates $x^s_{ij}$, $y^s_{ij}$ and $\xi^s_{ij}$, $\eta^s_{ij}$, respectively. (As above we assume that $x^s_{ij}$, $y^s_{ij}$, $\xi^s_{ij}$ and $\eta^s_{ij}$ are non-trivial. That is they are  not contained in the identity submatrix.)
Let us define the transition functions between two  superdomains $\mathcal U_I$ and $\mathcal U_J$ that correspond to $I = (I_s)$ and $J = (J_s)$, respectively, by the following formulas:
\begin{equation}\label{eq transition functions}
Z_{J_1} = Z_{I_1}C_{I_1J_1}^{-1}, \quad Z_{J_s} =
C_{I_{s-1}J_{s-1}}Z_{I_s}C_{I_sJ_s}^{-1},\quad s\ge 2. 
\end{equation}
Here $C_{I_sJ_s}$ is an invertible submatrix in $Z_{I_s}$ that coinsists of the
lines with numbers $i\in J_{s\bar 0}$ and $k_{s-1} + i,$ where $i\in J_{s\bar 1}$. In other words, we choose the matrix $C_{I_sJ_s}$ in such a way that $Z_{J_s}$ contains the identity submatrix $E_{k_s+l_s}$ in lines with numbers $i\in J_{s\bar 0}$ and $k_{s-1} + i,$ where $i\in J_{s\bar 1}$.
These charts and transition functions define a supermanifold that we denote by
$\mathbf F_{k|l}^{m|n}$. This supermanifold we will call the {\it flag supermanifold} of type $k|l$. In case $r = 1$ this supermanifold is called the {\it super-Grass\-mannian} and is denoted by $\mathbf
{Gr}_{m|n,k|l}$.

Let $\mathcal M = (\mathcal M_0,\mathcal O_{\mathcal M})$ be a complex-analytic supermanifold. 
Denote by 
$\mathcal T = \mathcal Der\,(\mathcal O_{\mathcal M})$  the sheaf of vector fields on  $\mathcal M$. It is a sheaf of Lie superalgebras with respect to the following multiplication 
$$
[X,Y] = XY
-(-1)^{p(X)p(Y)}YX.
$$
The global sections of $\mathcal T$
are called {\it holomorphic vector fields} on $\mathcal M$. They form a complex Lie superalgebra
that we denote by $\mathfrak v(\mathcal M)$. This Lie superalgebra is finite dimensional in case when $\mathcal M_0$ is compact. The goal of this paper is to compute the Lie superalgebra $\mathfrak v(\mathcal M)$ when $\mathcal M$ is a flag supermanifold of type $k|l$ in $\mathbb C^{m|n}$.

As usual we denote by $\mathfrak {gl}_{m|n}(\mathbb C)$ the general Lie superalgebra of the superspace $\mathbb C^{m|n}$. It coinsists of the following matrices:
$$
\left(
  \begin{array}{cc}
    A & B \\
    C & D \\
  \end{array}
\right), \quad\text{where} \quad  A\in \mathfrak{gl}_m(\mathbb C)\quad \text{and} \quad B\in \mathfrak{gl}_n(\mathbb C).
$$
Denote by $\GL_{m|n}(\mathbb C)$ the Lie supergroup of the Lie superalgebra $\mathfrak {gl}_{m|n}(\mathbb C)$. (See \cite{ViLieSupergroup} for more information about Lie supergroups.) In \cite{Man} an action of  $\GL_{m|n}(\mathbb C)$ on the supermanifold
 $\mathbf {F}_{k|l}^{m|n}$ was defined. In our coordinates this action is given by the following formulas:
\begin{equation}\label{eq action of Q}
\begin{aligned} 
&(L,(Z_{I_1},\ldots,Z_{I_r})) \longmapsto
(\tilde Z_{J_1},\ldots,\tilde Z_{J_r}), \quad \text{where}\\
 &L\in \GL_{m|n}(\mathbb C),\quad \tilde Z_{J_1} =
LZ_{I_1}C_1^{-1},\quad \tilde Z_{J_s} = C_{s-1}Z_{I_s}C_s^{-1}.
\end{aligned} 
\end{equation}
Here $C_1$ is an invertible submatrix in $LZ_{I_1}$ that consists of the lines with numbers $i\in J_{1\bar 0}$ and $m+i$, where  $i\in J_{1\bar 1}$, and $C_s,$ where $s\ge 2$, is an invertible submatrix in $C_{s-1}Z_{I_s}$ that consists of the lines with numbers $i\in J_{s\bar 0}$ and $k_{s-1}+i$, where
$i\in J_{\bar 1s}$. This Lie supergroup action induces a Lie superalgebra homomorphism 
$$
\mu:\mathfrak {gl}_{m|n}(\mathbb C)\to\mathfrak v(\mathbf {F}_{k|l}^{m|n}).
$$ 
In case
$r=1$ in \cite[Lemma 1]{onigl}  it was proven that $\Ker\mu = \langle E_{m+n}\rangle$, where
$E_{m+n}$ is the identity matrix of size $m+n$. In general case $r>1$ we also have $\Ker\mu = \langle E_{m+n}\rangle$ and the proof is similar to \cite[Lemma 1]{onigl}.  We see that $\mu$ induces an injective homomorphism of Lie superalgebras 
$$
\bar{\mu}:\mathfrak {gl}_{m|n}(\mathbb C)/\langle E_{m+n}\rangle\to\mathfrak v(\mathbf {F}_{k|l}^{m|n}).
$$
We will show that with some exceptions this homomorphism is an isomorphism.

\subsection{Superbundles and projectable vector fields}

Recall that a {\it morphism of complex-analytic} supermanifolds
$\mathcal M$ to  $\mathcal N$ is a pair $f = (f_0,f^*)$, where $f_0: \mathcal M_0\to \mathcal N_0$ is a holomorphic map and $f^*: \mathcal O_{\mathcal N}\to (f_0)_*(\mathcal
O_{\mathcal M})$ is a homomorphism of sheaves of superalgebras. 

\medskip
\de\label{de superbundle} A 
{\it superbundle} is a set $(\mathcal M, \mathcal B,p,\mathcal S)$, where $\mathcal S$ is fiber,  $\mathcal B$ is base space, $\mathcal M$ is total space and $p=(p_0,p^*):
\mathcal M\rightarrow \mathcal B$ is projection,   such that there exists an open covering 
 $\{U_i\}$ on $\mathcal B_0$, isomorphisms $\psi_i: (p_0^{-1}(U_i),\mathcal O_{\mathcal M})\rightarrow
(U_i,\mathcal {\mathcal O}_{\mathcal B})\times \mathcal S$ and the following diagram is commutative:
\begin{center}
\begin{tikzpicture}
\matrix (m) [matrix of math nodes,row sep=3em,column sep=4em,minimum width=2em]
 {(p_0^{-1}(U_i),\mathcal O_{\mathcal M}) 
 	& (U_i,\mathcal {\mathcal O}_{\mathcal B})\times \mathcal S \\
	(U_i,\mathcal {\mathcal O}_{\mathcal B})
	 & (U_i,
	 \mathcal {\mathcal O}_{\mathcal B})
	  \\};
\path[-stealth]
(m-1-1) edge node [left] {$p$} (m-2-1)
edge node [below] {$\psi_i$} (m-1-2)
(m-2-1.east|-m-2-2) edge node [below] {$id$}  (m-2-2)
(m-1-2) edge node [right] {$pr$} (m-2-2);
\end{tikzpicture}
\end{center}
where $pr$ is the natural projection.
\medskip

Usually we will denote a superbundle $(\mathcal M, \mathcal B,p,\mathcal S)$ just by $\mathcal M$. 

\medskip

\noindent {\bf Remark.} From the form of transition functions (\ref{eq transition functions}) it follows that for $r > 1$ the supermanifold $\mathbf {F}_{k|l}^{m|n}$ is a superbundle with base $\mathbf {Gr}_{m|n,k_1|l_1}$
and fiber $\mathbf {F}_{k'|l'}^{k_1|l_1}$, where $k'=(k_2,\ldots,k_r)$ and $l'=(l_2,\ldots,l_r)$. In local coordinates introduced above the projection $p$ is given by 
$$
(Z_1,Z_2,\ldots Z_n) \longmapsto (Z_1).
$$
Moreover, Formulas (\ref{eq action of Q}) tell us that the projection $p$ is equivariant with respect to the action of the supergroup  $\GL_{m|n}(\mathbb C)$ on $\mathbf {F}_{k|l}^{m|n}$ and $\mathbf {Gr}_{m|n,k_1|l_1}$.
\medskip

Let $p=(p_0,p^*): \mathcal M\to \mathcal N$ be a morphism of supermanifolds.

\medskip
\de \label{de projectable vector fields} A vector field $v\in\mathfrak
v(\mathcal M)$ is called {\it projectable} with respect to $p$,
if there exists a vector field $v_1\in\mathfrak v(\mathcal N)$
such that 
$$
p^*(v_1(f))=v(p^*(f))\quad \text{for all} \quad f\in \mathcal O_{\mathcal N}.
$$
In this case we say that {\it $v$ is projected to} $v_1$.

\medskip

Projectable vector fields form a super Lie subalgebra  $\overline{\mathfrak
v}(\mathcal M)$ in $\mathfrak v(\mathcal M)$.
In case if $p$ is a projection of a superbundle, the homomorphism
$p^*: \mathcal O_{\mathcal N}\to p_*(\mathcal O_{\mathcal M})$ is injective. Hence, any projectable vector field $v$ is projected into unique vector field
 $v_1 = \mathcal{P}(v)$ and we have the following map 
 $$
 \mathcal{P}:
\overline{\mathfrak v}(\mathcal M)\to \mathfrak v(\mathcal N),\quad v \mapsto v_1.
$$
It is a homomorphism of Lie superalgebras. A vector field $v\in\mathfrak v(\mathcal M)$ is called {\it vertical}, if $\mathcal{P}(v)=0$.
Vertical vector fields form an ideal  $\Ker\mathcal{P}$ in
Ëè $\overline{\mathfrak v}(\mathcal M)$.

We will need the following proposition proved in \cite{Bash}.

\medskip

\prop\label{prop Bash} {\sl Let  
	$p: \mathcal M 	\to \mathcal B$ be the projection of a superbundle with fiber $\mathcal S$. Assume that $\mathcal O_{\mathcal S}(\mathcal S_0) = \mathbb C$, i.e. any global holomorphic function is constant.   Then any holomorphic vector field from $\mathfrak v(\mathcal M)$ is projectable with respect to $p$ and we have a homomorphism of Lie superalgebras
	$	\nu : \mathfrak v(\mathcal M)\to \mathfrak v(\mathcal B).$
	}

\medskip

Let $p: \mathcal M\to
\mathcal B$ be a superbundle with fiber $\mathcal S$. We define the sheaf $\mathcal W$ on
$\mathcal B_0$ in the following way. We asign to any open set $U\subset \mathcal B_0$
 the set of all vertical vector fields on the supermanifold $(p_0^{-1}(U),\mathcal
O_{\mathcal M})$. In \cite{Vi} the following proposition was proven.

\medskip
\prop\label{prop W is localy free}  {\sl Assume that $\mathcal S_0$ is compact. Then $\mathcal W$ is a localy free sheaf of $\mathcal O_{\mathcal B}$-modules and
 $\dim\mathcal
W= \dim\mathfrak v(\mathcal S)$.}
\medskip

Clearly, the Lie algebra $\mathcal W(\mathcal B_0)$
coincides with the ideal of all vertical vector fields in $\mathfrak v(\mathcal M)$. Let us describe the corresponding to $\mathcal W$ graded sheaf.  Consider the following filtration in ${\mathcal O}_{\mathcal B}$
$$
{\mathcal O}_{\mathcal B} = \mathcal J^0 \supset \mathcal J^1 \supset \mathcal
J^2\dots
$$ 
where $\mathcal J$ is the sheaf of ideals in ${\mathcal O}_{\mathcal B}$ generated by odd elements. We have the corresponding graded sheaf of superalgebras
$$
\tilde{\mathcal O}_{\mathcal B}=\bigoplus_{p\ge 0}(\tilde{\mathcal O}_{\mathcal B})_{p},\quad \text{where}\quad (\tilde{\mathcal O}_{\mathcal B})_{p}=\mathcal J^p/\mathcal J^{p+1}.
$$
Putting $\mathcal W_{(p)}=
\mathcal J^p\mathcal W$ we get 
the following filtration in $\mathcal W$:
\begin{equation}\label{eq filtrations in W}
\mathcal W=\mathcal W_{(0)}
\supset\mathcal W_{(1)}\supset\dots .
\end{equation}
We define the $\mathbb Z$-graded 
sheaf of 
$\mathcal F_{\mathcal B_0}$-modules 
by
\begin{equation}\label{eq tilda W def}
\tilde {\mathcal W}=\bigoplus_{p\ge 0}\tilde{\mathcal W}_{p}, \quad
\text{where}\quad
\tilde{\mathcal W}_{p}=\mathcal W_{(p)}/\mathcal
W_{(p+1)}.
\end{equation}
Here $\mathcal F_{\mathcal B_0}$ is the structure sheaf of the underlying space $\mathcal B_0$. 
The ${\mathbb Z}_2$-grading in $\mathcal W_{(p)}$ induces the
${\mathbb Z}_2$-grading in $\tilde{\mathcal W}_{p}$.  Using Proposition \ref{prop W is localy free} we get the following result.

\medskip
\prop\label{prop tilde W_0 is localy free}  {\sl Assume that $\mathcal S_0$ is compact. Then $\tilde{\mathcal W}_0$ is a locally free sheaf of $\mathcal F_{\mathcal B_0}$-modules. Any fiber of the corresponding vector bundle is isomorphic to $\mathfrak v(\mathcal S)$.
	
}
\medskip

\subsection{The Borel-Weyl-Bott Theorem}

To calculate the Lie superalgebra of vector fields we will use the Borel-Weyl-Bott Theorem, see for example \cite{ADima} for details. This theorem permits to compute cohomology with values in a holomorphic homogeneous bundle over a flag manifold. For completeness we formulate it here adapting to our notations and agreements.  

Let  $G= \GL_{m}(\mathbb C)\times \GL_{n}(\mathbb C)$ be the underlying space of $\GL_{m|n}(\mathbb C)$,  $P$ be a parabolic subgroup in $G$ and $R$ be the reductive part of $P$. Assume that $\mathbf E_{\varphi}\to G/P$ is the homogeneous vector bundle corresponding to a representation  $\varphi$ of
$P$ in $E=(\mathbf E_{\varphi})_{P}$. Denote by $\mathcal E_{\varphi}$ the sheaf of holomorphic section of this vector bundle. 
In the Lie superalgebra  $\mathfrak{gl}_{m|n}(\mathbb C)_{\bar 0}\simeq \mathfrak{gl}_m(\mathbb C)\oplus \mathfrak{gl}_m(\mathbb C)$
we fix the Cartan subalgebra $
\mathfrak t=\mathfrak t_{1}\oplus \mathfrak t_{2},$
where
$$
\mathfrak t_{1}=
\{\operatorname{diag}(\mu_1,\dots,\mu_m)\}\quad\text{and} \quad \mathfrak t_{2}=
\{\operatorname{diag}(\lambda_1,\dots,\lambda_n)\},
$$
 the following system of positive roots:
 $$
 \Delta^+=\Delta^+_{1} \cup \Delta^+_{2},
 $$
 where
$$
\Delta^+_{1}=\{\mu_i-\mu_j, \,\,i<j\} \quad\text{and} \quad \Delta^+_{2}=\{\lambda_p-\lambda_q, \,\,p<q\},
$$
and the following system of simple roots $\Phi=\Phi_{1}\cup \Phi_{2}$, where
$$
\begin{array}{c}
\Phi_{1}= \{\alpha_1,..., \alpha_{n-1}\}, \,\,\,
\alpha_i=\mu_i-\mu_{i+1},\quad \Phi_{2}= \{\beta_1,..., \beta_{n-1}\}, \,\,\,
\beta_p=\lambda_p-\lambda_{p+1}.
\end{array}
$$
Denote by $\mathfrak t^*(\mathbb R)$
a real subspace in $\mathfrak t^*$
spaned by $\mu_j$ and $\lambda_i$. Consider the scalar product $( \,,\, )$ in $\mathfrak t^*(\mathbb R)$ such that the vectors  $\mu_j,\lambda_i$ form an orthonormal basis. An element $\gamma\in \mathfrak t^*(\mathbb R)$ is called {\it dominant} if $(\gamma, \alpha)\ge 0$ for all $\alpha \in \Delta^+$.

\medskip

\t [Borel-Weyl-Bott]. \label{teor borel} {\sl Assume that the representation	$\varphi: P\to \GL(E)$ is completely reducible and $\lambda_1,..., \lambda_s$
	are highest weights of $\varphi|R$. Then the $G$-module $H^0(G/P,\mathcal E_{\varphi})$ is isomorphic to the sum of irruducible $G$-modules with highest weights $\lambda_{i_1},..., \lambda_{i_t}$, where 
	$\lambda_{i_a}$ are dominant highest weights.

}

\medskip

\subsection{Holomorphic functions on flag supermanifolds}

Holomorphic functions on homogeneous supermanifolds and in particular on flag supermanifolds were studied in \cite{Viholom}. It is well-known that any holomorphic function on a connected compact complex
manifold is constant. This statement is false for a supermanifold with a
connected compact underlying space. However in case of flag supermanifolds the following theorem holds true: 

\medskip
\t\cite{Viholom}\label{teor constant functions} {\sl Consider the flag supermanifold $\mathcal M=\mathbf{F}_{k|l}^{m|n}$.  Assume that 
	\begin{equation}\label{eq_condition no functions}
	\begin{split}
	 (k|l)&\ne 
	(m,\ldots, m, k_{s+2},\ldots, k_r)| (l_1,\ldots, l_s,0,\ldots,0),\\
	(k|l)&\ne (k_1,\ldots, k_s,0,\ldots,0) |(n,\ldots, n, l_{s+2},\ldots,
	l_r),
\end{split}
\end{equation}
for any $s\geq 0$. Then $\mathcal O_{\mathcal M}(\mathcal M_0) = 	\mathbb{C}.$ In other words under conditions (\ref{eq_condition no functions}) any holomorphic function on $\mathbf{F}_{k|l}^{m|n}$ is constant.

 Otherwise
	$$
	\mathbf{F}_{k|l}^{m|n} \simeq
	(\pt,\bigwedge(mn))\times (\mathbf{F}^m_{k}\times \mathbf{F}^n_{l})
	$$
	and  $\mathcal O_{\mathcal M}(\mathcal M_0) = \bigwedge(mn)$, where $\bigwedge(mn)$ is the Grassmann algebra with $mn$ generators. 
}

\medskip

\section{Vector fields on flag supermanifolds}

\subsection{Vector fields on super-Grassmannians}

In previous sections we have seen that $\textbf{Gr}_{m|n,k|l}$ is a $\GL_{m|n}(\mathbb{C})$-homogene\-ous superspace. The action of $\GL_{m|n}(\mathbb{C})$ on $\textbf{Gr}_{m|n,k|l}$ is given by Formulas (\ref{eq action of Q}). This action induces the Lie algebra homomorphism 
$$
\mu: \mathfrak {gl}_{m|n}(\mathbb C)\to
\mathfrak{v}(\textbf{Gr}_{m|n,k|l}).
$$
The kernel of this homomorphism is eqaul to $\langle E_{m+n}\rangle$, \cite[Lemma 1]{onigl}.
Further we will use the following notation: 
$$
\mathfrak{pgl}_{m|n}(\mathbb C):=\mathfrak {gl}_{m|n}(\mathbb C)/
\langle E_{m+n}\rangle.
$$
The Lie superalgebra of holomorphic vector fields on  super-Grassmannian $\mathbf{Gr}_{m|n,k|l}$ was computed in \cite{buneg,onigl,onich,serov}.

\medskip
\t\label{teor super=Grassmannians} {\sl 
			The homomorphism  $\mu:\mathfrak {gl}_{m|n}(\mathbb C)\to
			\mathfrak{v}(\textbf{Gr}_{m|n,k|l})$
	is almost always surjective and  
	$$
	\mathfrak v(\textbf{Gr}_{m|n,k|l})\simeq \mathfrak{pgl}_{m|n}(\mathbb C).
	$$
	The exeptional cases are the following. 
	
	\begin{description}
		\item[1.1] For the super-Grassmannian $\mathbf{Gr}_{2|2,1|1}$ we have 
		$$
		\mathfrak v(\mathbf{Gr}_{2|2,1|1})\simeq
		\mathfrak {psl}_{2|2}(\mathbb C)+\!\!\!\!\!\!\supset\mathfrak {sl}_{2}(\mathbb
		C),
		$$
		where $\mathfrak {psl}_{2|2}(\mathbb C)=\mathfrak {sl}_{2|2}(\mathbb
		C)/{<E_4>}$.
		
		\item[1.2] For $\mathbf{Gr}_{1|n,0|n-1}\simeq\mathbf{Gr}_{n|1,n-1|0}\simeq \mathbf{Gr}_{n|1,1|1} \simeq\mathbf{Gr}_{1|n,1|1}$, $n>2$,   we have
		$$
		\mathfrak v(\mathbf{Gr}_{1|n,0|n-1})\simeq W_n =
		\operatorname{Der}\bigwedge(\zeta_1,\dots,\zeta_n).
		$$
		
		\item[1.3] In the degenerate case $\mathbf{Gr}_{m|n,0|n}\simeq
		\mathbf{Gr}_{m|n,m|0}$ we have
		$$
		\mathfrak v(\mathbf {F}^{m|n}_{0|n})\simeq W_{mn}=
		\operatorname{Der}\bigwedge (\zeta_1,\dots,\zeta_{mn}).
		$$
		
		\item[1.4] For $\mathbf{Gr}_{2|2,0|1}\simeq \mathbf{Gr}_{2|2,1|0}
		\simeq\mathbf{Gr}_{2|2,1|2}\simeq \mathbf{Gr}_{2|2,2|1}$ we have
		$$
		\mathfrak v(\mathbf{Gr}_{2|2,0|1})\simeq \tilde{\mathbf
			H}_4+\!\!\!\!\!\!\supset\langle z\rangle,
		$$
		where $\ad z$ acts on the Lie superalgebra of Cartan type $\tilde{\mathbf H}_4$ as the grading operator. 
	\end{description}

}

\medskip

In case
$$
 0<k<m\quad  \text{and}\quad  0<l<n, \quad (m|n,k|l)\ne (2|2,1|1),
$$
  the Lie superalgebra of vector fields was computed in \cite{onigl}. Results  $1.1$ and $1.2$ of Theorem \ref{teor super=Grassmannians}  were obtained in \cite{buneg} (see also \cite{onigl} for an explicit description of the Lie superalgebra) and \cite{onich,serov}, respectively. Result   $1.3$ of Theorem \ref{teor super=Grassmannians}  is obvious.
Result $1.4$ of Theorem \ref{teor super=Grassmannians} follows from arguments in \cite{serov}, Proof of Theorem $2.6$. Note that in the statement of Theorem $2.6$ in \cite{serov} and also in \cite[Theorem 7]{onigl}  the Lie superalgebra of vector fields in case $1.4$ was pointed incorrectly.

We will need an explicit description of the Lie superalgebra of holomorphic vector fields on  $\mathbf{Gr}_{2|2,1|1}$, Case $1.1$ of Theorem \ref{teor super=Grassmannians},  in the following local chart  
$$
\left(
\begin{array}{cc}
x & \xi \\
1 & 0  \\
\eta & y \\
0 & 1 \\
\end{array}
\right).
$$
The image of $\mathfrak{gl}_{2|2}(\mathbb C)_{\bar 0}$ with respect to the homomorphism $\mu$ from Theorem \ref{teor super=Grassmannians} is given by:
\begin{equation}\label{eq vector fields GR_2|2,1|1 even}
\begin{split}
&\mu(E_{11}) = x\frac{\partial}{\partial x} + \xi\frac{\partial}{\partial \xi}, \quad  \mu(E_{12}) = \frac{\partial}{\partial x}, \quad \mu(E_{22}) = - x\frac{\partial}{\partial x}- \eta \frac{\partial}{\partial \eta},\\
 &\mu(E_{21}) = -x^2\frac{\partial}{\partial x}-x\eta \frac{\partial}{\partial \eta} - x\xi \frac{\partial}{\partial \xi} + \xi\eta \frac{\partial}{\partial y} ,\quad \mu(E_{34}) = \frac{\partial}{\partial y}, \\
 & \mu(E_{43}) = -y^2\frac{\partial}{\partial y}-y\xi \frac{\partial}{\partial \xi} - y\eta \frac{\partial}{\partial \eta} - \xi\eta \frac{\partial}{\partial x} ,\quad \mu(E_{33}) = y\frac{\partial}{\partial y} + \eta\frac{\partial}{\partial \eta},\\
 & \mu(E_{44}) = - y\frac{\partial}{\partial y}- \xi \frac{\partial}{\partial \xi}.
\end{split}
\end{equation}
The image of $\mathfrak{gl}_{2|2}(\mathbb C)_{\bar 1}$ with respect to the homomorphism $\mu$ from Theorem \ref{teor super=Grassmannians} is given by:
\begin{equation}\label{eq vector fields GR_2|2,1|1 odd}
\begin{split}
&\mu(E_{14}) = \frac{\partial}{\partial \xi},\quad  \mu(E_{32}) = \frac{\partial}{\partial \eta}, \quad  \mu(E_{13}) = \eta\frac{\partial}{\partial x} + y\frac{\partial}{\partial \xi},\\
&  \mu(E_{31}) = \xi\frac{\partial}{\partial y} + x\frac{\partial}{\partial \eta}, \quad \mu(E_{23}) = - x\eta \frac{\partial}{\partial x} - xy \frac{\partial}{\partial \xi} + y\eta \frac{\partial}{\partial y},\\
& \mu(E_{41}) = - y\xi \frac{\partial}{\partial y} - xy \frac{\partial}{\partial \eta} + x\xi \frac{\partial}{\partial x},\quad 
\mu(E_{24}) = -x\frac{\partial}{\partial \xi} + \eta\frac{\partial}{\partial y},\\
& \mu(E_{42}) = -y\frac{\partial}{\partial \eta} + \xi\frac{\partial}{\partial x}.
\end{split}
\end{equation}
Additional holomorphic on $\mathbf{Gr}_{2|2,1|1}$ vector fields are 
\begin{equation}\label{eq vector fields GR_2|2,1|1 additional}
\begin{split}
&\eta\frac{\partial}{\partial \xi}, \quad \xi\frac{\partial}{\partial \eta}.
\end{split}
\end{equation}
A direct computation shows that 
\begin{equation}\label{eq Gr_2|2,1|1 decomposition of Gl_0 mod}
\mathfrak v(\mathbf{Gr}_{2|2,1|1}) \simeq \mathfrak{pgl}_{2|2}(\mathbb C)_{\bar 0}\oplus \mathfrak{pgl}_{2|2}(\mathbb C)_{\bar 1}\oplus \langle \eta\frac{\partial}{\partial \xi}, \,\, \xi\frac{\partial}{\partial \eta} \rangle
\end{equation}
as $\mathfrak{gl}_{2|2}(\mathbb C)_{\bar 0}$-modules.

Let us give an explicit description of the Lie superalgebra of holomorphic vector fields on  $\mathbf{Gr}_{2|2,1|2}$, Case $1.4$ of Theorem \ref{teor super=Grassmannians}  in the following local chart  
\begin{equation}\label{eq chart on F_2|2,1|2}
\left(
\begin{array}{ccc}
x & \xi_1 & \xi_2 \\
1 & 0 & 0 \\
0 & 1 & 0 \\
0 & 0 & 1 \\
\end{array}
\right).
\end{equation}
The definition of the Lie superalgebra $\tilde{\mathbf H}_4$ can be found  in \cite{Kac}. For completeness we remind it here. We have $\tilde{\mathbf H}_4\subset \operatorname{Der}\bigwedge(\theta_1,\dots,\theta_4)$ and $\tilde{\mathbf H}_4$ consists of all elements in the form:
$$
D_f=\sum_{i=1}^4\frac{\partial f}{\partial \theta_i} \frac{\partial }{\partial \theta_i}, \quad f\in \bigwedge(\theta_1,\dots,\theta_4), \quad f(0)=0.
$$
The Lie superalgebra $\tilde{\mathbf H}_4$ is $\mathbb Z$-graded and in chosen chart the image of an injective homomorphism $\tilde{\mathbf H}_4 \to 	\mathfrak v(\mathbf{Gr}_{2|2,1|2})$ is given by the following vector fields:
\begin{equation}\label{eq H_4 =}
\begin{split}
&(\tilde{\mathbf H}_4)_{-1}=  \Bigl\langle\frac{\partial}{\partial
	\xi_1}, \quad \frac{\partial}{\partial \xi_2}, \quad
x\frac{\partial}{\partial \xi_1}, \quad x\frac{\partial}{\partial \xi_2}\Bigr\rangle; \\
& (\tilde{\mathbf H}_4)_{0}=  \Bigl\langle \frac{\partial}{\partial x}, \quad x\frac{\partial}{\partial x} + \xi_1\frac{\partial}{\partial
	\xi_1}, \quad x\frac{\partial}{\partial x} + \xi_2\frac{\partial}{\partial
	\xi_2}, \quad \xi_1\frac{\partial}{\partial
	\xi_2}, \\
&\quad  \quad \quad \quad \quad\xi_2\frac{\partial}{\partial
	\xi_1}, \quad
 x\xi_1\frac{\partial}{\partial
	\xi_1}+x\xi_2\frac{\partial}{\partial
	\xi_2}+x^2\frac{\partial}{\partial x} \Bigr\rangle;\\
&(\tilde{\mathbf H}_4)_{1}=  \Bigl\langle \xi_1\frac{\partial}{\partial x}, \quad
\xi_2\frac{\partial}{\partial x}, \quad x\xi_1\frac{\partial}{\partial x}+ \xi_1\xi_2\frac{\partial}{\partial \xi_2},
\quad x\xi_2\frac{\partial}{\partial x}-
\xi_1\xi_2\frac{\partial}{\partial \xi_1} \Bigr\rangle; \\
&(\tilde{\mathbf H}_4)_{2}=  \Bigl\langle \theta = \xi_1\xi_2\frac{\partial}{\partial x}\Bigr\rangle; \\
\end{split}
\end{equation}
The $\mathbb Z$-graded operator mentioned in Theorem \ref{teor super=Grassmannians} is given by: 
$$
z= \xi_1\frac{\partial}{\partial \xi_1} + \xi_2\frac{\partial}{\partial \xi_2}.
$$

\medskip

We will call the super-Grassmannians from 1.1-1.4 of Theorem  \ref{teor super=Grassmannians} {\it exceptional}. Note that the super-Grassmannian $\mathbf{Gr}_{0|n,0|l} \simeq \mathbf{Gr}_{n|0,l|0}$ is just usual Grassmannians isomorphic to $\mathbf{Gr}_{n,l}$. It is well-known that 
$$
\mathfrak v(\mathbf{Gr}_{n,l})\simeq\mathfrak{pgl}_{n}(\mathbf C),
$$
see \cite{ADima} for details.

\subsection{Vector fields on flag supermanifolds. Main case}

Assume that $r>1$. From now on we use the following notations: 
$$
\mathcal M=\mathbf{F}^{m|n}_{k|l}, \quad \mathcal B= 
\mathbf{Gr}_{m|n,k_1|l_1} \quad \text{and} \quad  \mathcal S=\mathbf{F}_{k'|l'}^{k_1|l_1},
$$
where $k'= (k_2,\ldots, k_r)$ and $l'= (l_2,\ldots, l_r)$.  If $\mathcal O_{\mathcal S}(\mathcal S_0)=\mathbb C$, then by Proposition \ref{prop Bash} the projection of the superbundle
$\mathcal M\to \mathcal B$ determines the homomorphism of Lie superalgebras 
$$
\mathcal{P}:\mathfrak v(\mathcal M)\to\mathfrak
v(\mathcal B).
$$
This projection is $\GL_{m|n}(\mathbb C)$-equivariant. Hence for the natural Lie superalgebra homomorphisms $\mu: \mathfrak{gl}_{m|n}(\mathbb C)\to\mathfrak
v(\mathcal M)$ and $\mu_{\mathcal B}: \mathfrak{gl}_{m|n}(\mathbb C)\to\mathfrak v(\mathcal B)$ we have 
$$
\mu_{\mathcal B} = \mathcal{P}\circ\mu.
$$
By Theorem \ref{teor super=Grassmannians}, the homomorphisms $\mu_{\mathcal B}$ and hence the homomorphism
$\mathcal{P}$ is almost always surjective. We will prove that $\mathcal{P}$ is injective. Hence, 
\begin{equation}\label{eq connection mu and mu_B}
\mu = \mathcal{P}^{-1}\circ\mu_{\mathcal B}
\end{equation}
is surjective and  
$$
\mathfrak v(\mathcal M)\simeq\mathfrak{gl}_{m|n}(\mathbb C)/ \langle
E_{m+n}\rangle .
$$

 In previous section 
we constructed a locally free sheaf $\tilde{\mathcal
W}$ on $\mathcal B_0$. The sheaf $\mathcal W$ possesses the natural action of the Lie group $G = \GL_m(\mathbb C)\times
\GL_n(\mathbb C)$, because $G$ is the underlying space of $\GL_{m|n}(\mathbb C)$. This action preserves the filtration (\ref{eq filtrations in W}) and induces the action 
in the sheaf $\tilde{\mathcal W}$. Hence the vector bundle $\mathbf W_0\to \mathcal B_0$ corresponding to the localy free sheaf
$\tilde{\mathcal W}_0$ is homogeneous. 
Consider the local chart on the
super-Grassmannian $\mathcal B$ corresponding to 
\begin{equation}\label{eq I_1 of o in  B_0}
I_{1\bar 0} = \{m-k_1+1,\ldots,m\}\quad \text{and}\quad I_{1\bar 1} = \{n-l_1+1,\ldots,n\}.
\end{equation}
The coordinate matrix $Z_{I_1}$ in this case has the following form
\begin{equation}\label{eq loc chart on B}
Z_{I_1} =
\left(
\begin{array}{cc}
X_1&\Xi_1\\
E_{k_1}&0\\
\H_1&Y_1\\0&E_{l_1}\end{array} \right). 
\end{equation}
Denote by $o$ the point in $\mathcal B_0$ defined by the following equations:
$$
X_1=Y_1=\Xi_1=\H_1=0.
$$
Then $\mathcal B_0$ is naturally isomorphic to $G/H$, where $H$ is the stabilizer of $o$. An easy computation shows that $H$ contains all matrices in the following form:
\begin{equation}\label{eq H contais matrices}
\left(
\begin{array}{cccc} 
A_1&0&0 &0\\
C_1&B_1&0&0\\
0&0&A_2&0\\
0&0&C_2&B_2
\end{array} 
\right),
\end{equation} 
where 
$$
A_1\in \GL_{m-k_1}(\mathbb C),\,\, A_2\in \GL_{n-l_1}(\mathbb C),\,\, B_1\in \GL_{k_1}(\mathbb C)\,\, \text{and} \,\, B_2\in \GL_{l_1}(\mathbb C).
$$
The reductive part $R$ of $H$ is given by the following equations 
$$
C_i=0,\quad i=1,2.
$$

 Let us compute the representation $\psi$ of $H$ in the fiber $(\mathbf W_0)_o$ of $\mathbf W_0$ over the point $o$. We identify $(\mathbf W_0)_o$ with the Lie superalgebra of holomorphic vector fields $\mathfrak v(\mathcal S)$, see Proposition \ref{prop tilde W_0 is localy free}. 
Let us choose an atlas on $\mathcal M$  defined by $I_{1}=(I_{1\bar 0}, I_{s\bar 1})$, see (\ref{eq I_1 of o in  B_0}), and by certain $I_{s},\; s = 2,\ldots,r$. In notations (\ref{eq loc chart on B}) and (\ref{eq H contais matrices}) the group $H$ acts in the chart defined by $Z_{I_1}$ in the following way:
$$
\left(
\begin{array}{cccc} 
A_1&0&0 &0\\
C_1&B_1&0&0\\
0&0&A_2&0\\
0&0&C_2&B_2
\end{array} 
\right)
\left(
\begin{array}{cc}
X_1&\Xi_1\\
E_{k_1}&0\\
\H_1&Y_1\\
0&E_{l_1}
\end{array} \right) 
= 
\left(
\begin{array}{cc} 
A_1 X_1&A_1\Xi_1\\
 C_1 X_1 + B_1 & C_1\Xi_1\\
A_2\H_1 & A_2 Y_1\\ 
C_2\H_1 & C_2 Y_1 + B_2 
\end{array} \right).
$$
Hence, for $Z_{I_2}$ we have
\begin{equation}\label{eq action over o}
\begin{split} 
&\left(
\begin{array}{cc} 
C_1 X_1 + B_1 & C_1\Xi_1\\
C_2\H_1 & C_2 Y_1 + B_2
 \end{array} \right)
\left(
\begin{array}{cc}  
X_2 & \Xi_2 \\
 \H_2 & Y_2
\end{array} 
\right) =\\
&\quad \quad = \left(
\begin{array}{cc} 
(C_1 X_1 + B_1)X_2+ C_1\Xi_1\H_2  & (C_1 X_1 + B_1)\Xi_2+ C_1\Xi_1Y_2
\\
C_2\H_1X_2+ (C_2 Y_1 + B_2)\H_2 & C_2\H_1\Xi_2+ (C_2 Y_1 + B_2)Y_2
\end{array} \right).
\end{split}
\end{equation}
Note that the local coordinates of $Z_{I_s},\; s \ge 2$, can be interpreted as local coordinates on the fiber $\mathcal S$ of the superbundle $\mathcal M$. To obtain the action of $H$ in the fiber $(\bf W_0)_o$ in these coordinates we put  
$$
X_1=Y_1 = 0 \quad \text{and} \quad \Xi_1=\H_1=0
$$
in (\ref{eq action over o}) and  modify
$Z_{I_s},\; s \ge 3$, accordingly.
We see that the nilradical of $H$ and the subgroup
$\GL_{m-k_1}(\mathbb C)\times \GL_{n-l_1}(\mathbb C)$ in $R$ act trivially on
$\mathcal S$ and that  the subgroup $\GL_{k_1}(\mathbb C)\times \GL_{l_1}(\mathbb C)\subset
R$ acts in the natural way. In other words the action of $H$ in $\mathcal S$ over $o$ is given by the following formulas:
\begin{equation}
\label{eq Z_I_2_preobrazuetsja}
\begin{split}
\quad\left(
\begin{array}{cc}  B_1 & 0\\ 0 &  B_2\end{array}
\right) \left(
\begin{array}{cc}  X_2 & \Xi_2 \\ \H_2 & Y_2\end{array} \right) =
\left(
\begin{array}{cc} B_1X_2  &
B_1\Xi_2\\
B_2\H_2  & B_2 Y_2
\end{array} \right).
\end{split}
\end{equation}
This means that $H$ acts as the underlying space of the Lie supergroup $\GL_{k_1|l_1}(\mathbb C)$ on the flag supermanifold $\mathcal S$, see (\ref{eq action of Q}). 
Furthermore  assume that 
$$
\mathfrak v(\mathcal S)\simeq \mathfrak{gl}_{k_1|l_1}(\mathbb C)/ \langle
E_{k_1+l_1}\rangle = \left\lbrace\left(
\begin{array}{cc}
Z_1 & T_1  \\
T_2 & Z_2  \\
\end{array}
\right)+<E_{k_1+l_1}>\right\rbrace,
$$
where $Z_1\in \mathfrak{gl}_{k_1}(\mathbb C)$ and $Z_2\in \mathfrak{gl}_{l_1}(\mathbb C)$. Then the induced action of the Lie group $\GL_{k_1}(\mathbb C)\times \GL_{l_1}(\mathbb C)$ on  $(\bf W_0)_o = \mathfrak v(\mathcal S)$ coinsides with the adjoint action of the underlying Lie group of the Lie supergroup $\GL_{k_1|l_1}(\mathbb
C)$. More precisely, we have
\begin{equation}\label{eq submodules of pgl}
\begin{split}
\left(
\begin{array}{cc}
B_1 & 0 \\
0 & B_2
\end{array}
\right)& \left( \left(
\begin{array}{cc}
Z_1 & T_1  \\
T_2 & Z_2  \\
\end{array}
\right)+<E_{k_1+l_1}>\right)
\left(
\begin{array}{cc}
B_1^{-1} & 0 \\
0 & B_2^{-1}
\end{array}
\right)=\\
&
\left(
\begin{array}{cc}
B_1Z_1B_1^{-1} & B_1T_1B_2^{-1} \\
B_2T_2B_1^{-1} & B_2Z_2B_2^{-1}
\end{array}
\right)+<E_{k_1+l_1}>,
\end{split}
\end{equation}
where $B_1\in \GL_{k_1}(\mathbb C)$ and $B_2\in \GL_{l_1}(\mathbb C)$.

Denote by $\Ad_{k_1}$ and $\Ad_{l_1}$ the adjoint representations of $\GL_{k_1}(\mathbb C)$ and $\GL_{l_1}(\mathbb C)$ on $\mathfrak {sl}_{k_1}(\mathbb C)$ and $\mathfrak {sl}_{l_1}(\mathbb C)$, respectively, and by $\rho_{k_1}$ and $\rho_{l_1}$ the standard representations of $\GL_{k_1}(\mathbb C)$ and $\GL_{l_1}(\mathbb C)$ in $\mathbb {C}^{k_1}$ and $\mathbb {C}^{l_1}$, respectively. We denote by
$1$ the one dimensional trivial representation of $\GL_{k_1}(\mathbb C)\times \GL_{l_1}(\mathbb C)$.
The following lemma follows from (\ref{eq submodules of pgl}).

\medskip

\lem\label{lem representation of Gl_k_1 in v(S)} 
 {\sl The representation $\psi$ of $H$ in the fiber
$(\mathbf W_0)_o = \mathfrak v(\mathcal S)$ is completely reducible. The nilradical of $H$ acts trivially in $(\mathbf W_0)_o$. If
$\mathfrak v(\mathcal S)\simeq \mathfrak{gl}_{k_1|l_1}(\mathbb C) / \langle
E_{k_1+l_1}\rangle $, 
then
\begin{equation}
\label{eq rep of H in v(S) general}
\psi|R =\left\{
\begin{array}{l}
\Ad_{k_1} + \Ad_{l_1} + \rho_{k_1}\otimes\rho_{l_1}^* +
\rho_{l_1}\otimes\rho_{k_1}^* + 1\,\,\, \text{for}\,\, k_1,l_1 > 0,\\
\Ad_{k_1}\,\,\, \text{for}\; k_1>0,\,\, l_1 = 0,\\
\Ad_{l_1}\,\,\, \text{for}\,\, k_1 = 0,\,\, l_1 > 0.\end{array}\right.
\end{equation}
}

\medskip

Further we will use the chart on $\mathbf {F}_{k|l}^{m|n}$
defined by 
$I_s= I_{s\bar 0}\cup I_{s\bar 1}$, where $I_{1\bar i}$ is as above, and  
$$
I_{s\bar 0}=\{k_{s-1}-k_s+1,\ldots,k_{s-1}\},\quad I_{s\bar 1}=\{l_{s-1}-l_s+1,\ldots,l_{s-1}\}
$$
for $s\ge 2$.
The coordinate matrices of this chart have the following form
$$
Z_{I_s} = \left(
\begin{array}{cc} 
X_s&\Xi_s\\ 
E_{k_s}&0\\
\H_s&Y_s\\
0&E_{l_s}
\end{array}\right),\quad s=1,\ldots k,
$$
where again the local coordinates are 
$$
X_s=(x^s_{ij}),\,\, Y_s=(y^s_{ij}), \,\, \Xi_s=(\xi^s_{ij})\,\,\, \text{and} \,\,\, \H_s=(\eta^s_{ij}).
$$
We denote this chart by $\mathcal U$ and the corresponding chart on $\mathcal B$ by $\mathcal U_{\mathcal B}$. In other words, $\mathcal U_{\mathcal B}$ is given by the coordinate matrix (\ref{eq loc chart on B}).

\medskip
\lem \label{lem some fundamental vector fields} {\sl The  vector fields $\frac{\partial}{\partial \xi^1_{ij}}$ and $\frac{\partial}{\partial \eta^1_{ij}}$ are fundamental. This is they are induced by the natural action of  $\GL_{m|n}(\mathbb C)$ on $\mathcal M$.
}
\medskip

\noindent{\it Proof.} Let us prove this statement for example for the vector field
$\frac{\partial}{\partial \xi^1_{11}}$. This vector field corresponds to the one-parameter subgroup $\operatorname{exp}(\tau
E_{1,a})$, where $a=m +n-l_1+1$ and $\tau$ is an odd parameter.
Indeed, the action of this subgroup is given by
$$
\left(\begin{array}{cc} 
X_1&\Xi_1\\ 
E_{k_1}&0\\
\H_1&Y_1\\
0&E_{l_1}
\end{array}\right)
\mapsto
\left(\begin{array}{cc}
 X_1&\tilde \Xi_1\\ 
E_{k_1}&0\\
\H_1&Y_1\\
0&E_{l_1}
\end{array}\right)\quad \text{and}\quad Z_{I_s}\mapsto Z_{I_s},\;
s\ge 2,
$$
where
$$
\tilde{\Xi}_1=\left(
\begin{array}{ccc}
\tau+\xi_{11}^1 & \ldots & \xi_{1l_1}^1 \\
\vdots & \ddots & \vdots \\
\xi_{m-k_1,1}^1 & \ldots & \xi_{m-k_1, l_1}^1 \\
\end{array}\right). \Box
$$

\medskip

Let us choose a basis $v_i$, where $i=1,\ldots \dim(\mathfrak v(\mathcal S))$,  in $\mathfrak v(\mathcal
S)$. Any holomorphic vertical vector field on $\mathcal M$ can be written uniquely in the form 
\begin{equation}\label{eq form of vertical vector fields}
w=\sum_q f_qv_q, 
\end{equation}
where $f_q$ are holomorphic functions on $\mathcal U$ depending only on coordinates from $Z_{I_1}$. We will need the following two lemmas:

\medskip
\lem \label{lem assume ker Pi ne 0} {\sl If $\Ker \mathcal{P}\ne \{0\}$, then $\dim \mathcal{W}_{(0)}(\mathcal B_0)> \dim  \mathcal{W}_{(1)}(\mathcal B_0).$	}

\medskip

Note that since $\mathcal B_0$ is compact, $\dim \mathcal{W}_{(i)}(\mathcal B_0) < \infty$ for all $i$.

\medskip

\noindent{\it Proof.} By definition we have the inclusion of sheaves $\mathcal{W}_{(1)} \hookrightarrow \mathcal{W}_{(0)}$ and hence we have the inclusion of the vector spaces of global sections 
$$
\mathcal{W}_{(1)}(\mathcal B_0) \hookrightarrow \mathcal{W}_{(0)}(\mathcal B_0).
$$
Therefore we need to show that there exists a vector field $v\in \mathcal{W}_{(0)}(\mathcal B_0)$ such that $v\notin \mathcal{W}_{(1)}(\mathcal B_0)$. Consider a vector field $w\in \mathcal{W}_{(1)}(\mathcal B_0)$ written in the form (\ref{eq form of vertical vector fields}). 
 Assume that there is a function $f_q$ that depends for example on odd coordinate $\xi_{ij}^1$. Then $w = \xi_{ij}^1 w'+ w''$, where $w'$
and $w''$ are local vertical vector fields and their coefficients (\ref{eq form of vertical vector fields}) do not depend on $\xi_{ij}^1$, and $w'\ne 0$. Using Lemma \ref{lem some fundamental vector fields} and the fact that $\operatorname{Ker}\mathcal{P}$ is an ideal in
$\mathfrak v(\mathcal M)$, we see that 
$$
w' =
[\frac{\partial}{\partial \xi^1_{ij}},w]\in
\operatorname{Ker}\mathcal{P}.
$$
In particular, $w'$ is a global vertical vector field. In this way we can exclude  all odd coordinates $\xi_{ij}^1$ and $\eta_{ij}^1$. Therefore there exists a vector field $v$ from $\Ker \mathcal{P}$ such that $v\in \mathcal{W}_{(0)}(\mathcal B_0)$ but $v\notin \mathcal{W}_{(1)}(\mathcal B_0)$.$\Box$

\medskip

\lem \label{lem sections of W_0} {\sl We have
\begin{equation}
\label{H0(F(gl))}
\widetilde{\mathcal W}_0(\mathcal B_0)\simeq \left\{
\begin{array}{ll}
\mathbb C, & 0<k_1<m,\,\,0<l_1<n ; \\
\mathfrak r_1\oplus\mathfrak r_2\oplus \mathbb C, & 1<k_1=m,\,\,0<l_1<n; \\
\mathfrak r_3\oplus\mathfrak r_4\oplus \mathbb C, & 0<k_1<m,\,\,1<l_1=n;\\
\mathfrak r_2\oplus \mathbb C, & 1=k_1=m,\,\,0<l_1<n; \\
\mathfrak r_3\oplus \mathbb C, & 0<k_1<m,\,\,1=l_1=n;\\
\{0\}, & 0<k_1<m,\,\,0=l_1\leq n,\,\,\text{or}\\
&0=k_1\leq m,\,\,0<l_1<n,\,\,\text{or}\\
&0=k_1< m,\,\,1=l_1\leq n,\,\,\text{or}\\
&1=k_1\leq m,\,\,0=l_1< n;\\
\mathfrak r_1, & 1<k_1=m,\,\,0=l_1<n;\\
\mathfrak r_4, & 0=k_1<m,\,\,1<l_1=n,\\
\end{array}
\right.
\end{equation}
where $\mathfrak r_1$, $\mathfrak r_2$, $\mathfrak r_3$, $\mathfrak
r_4$ are irreducible $\mathfrak {sl}_{m}(\mathbb C)\oplus \mathfrak
{sl}_{n}(\mathbb C)$-modules with the  highest weights
$\mu_{1}-\mu_{m}$, $\mu_{1}-\lambda_{n}$, $\lambda_{1}- \mu_{m}$ and $\lambda_{1}- \lambda_{n}$
respectively. The trivial $1$-dimensional module $\mathbb C$ corresponds to the highest weight $0$.	
	
}
\medskip

\noindent{\it Proof.} We compute the vector space of global sections of $\mathbf W_0$ using the Borel-Weyl-Bott Theorem \ref{teor borel}.
The representation $\psi$ of $H$ in $(\mathbf W_0)_o$ is described in Lemma \ref{lem representation of Gl_k_1 in v(S)}. From (\ref{eq rep of H in v(S) general}) it follows that the highest weights of $\psi$
have the form: 
\begin{itemize}
	\item  $\mu_{m-k_1+1}-\mu_{m}$, $\mu_{m-k_1+1}-\lambda_{n}$,
	$\lambda_{n-l_1+1}- \mu_{m}$, $\lambda_{n-l_1+1}- \lambda_{n}$, $0$ for  $k_1,\,l_1>1$;
	
	\item $\mu_{m}-\lambda_{n}$,
	$\lambda_{n-l_1+1}- \mu_{m}$, $\lambda_{n-l_1+1}- \lambda_{n}$, $0$ for  $k_1=1$, $l_1>1$;
	
	\item $\mu_{m-k_1+1}-\mu_{m}$, $\mu_{m-k_1+1}-\lambda_{n}$, $\lambda_{n}- \mu_{m}$, $0$ for  $k_1>1$, $l_1=1$;
	
	\item $\mu_{m}-\lambda_{n}$, $\lambda_{n}- \mu_{m}$, $0$ for  $k_1=1$, $l_1=1$;
	
	\item $\mu_{m-k_1+1}-\mu_{m}$ for $k_1>1$, $l_1=0$;
	\item $\lambda_{n-l_1+1}- \lambda_{n}$ for $k_1=0$, $l_1>1$.
\end{itemize}
(Note that for $k_1=1$, $l_1=0$ and $k_1=0$, $l_1=1$ the representation space of $\psi$ is trivial.)
Therefore the dominant highest weights of $\psi$
have the following form: 
\begin{itemize}
	\item $0$, if $0<k_1<m$ and $0<l_1<n$;
	
	\item $0$, $\mu_{1}-\mu_{m}$, $\mu_{1}-\lambda_{n}$, if $1<k_1=m$, $0< l_1<n$;
	\item $0$, $\mu_{1}-\lambda_{n}$, if $1=k_1=m$, $0< l_1<n$;

	\item $0$, $\lambda_{1}-\lambda_{n}$, $\lambda_{1}-\mu_{m}$, if $0<k_1<m$, $1<l_1=n$;
	\item $0$, $\lambda_{1}-\mu_{m}$, if $0<k_1<m$, $1= l_1=n$;
	
	\item $\mu_{1}-\mu_{m}$, if $1<k_1=m$, $0=l_1<n$;
	\item $\lambda_{1}- \lambda_{n}$, if $0=k_1<m$, $1<l_1=n$.
\end{itemize}
We have no dominant weights in the following cases: 
\begin{itemize}
	\item $0<k_1<m$, $0=l_1\leq n$;
\item $0=k_1\leq m$, $0<l_1<n$;
\item $0=k_1< m,\,\,1=l_1\leq n$;
\item  $1=k_1\leq m,\,\,0=l_1< n$.
\end{itemize}
By Borel-Weyl-Bott Theorem we get the result.$\Box$

\medskip

We are ready to prove the following theorem. 

\medskip
\t \label{teor Ker Pi = 0} {\sl Assume that  $r>1$. If 
	$$
	\mathcal O_{\mathcal S}(\mathcal S_0)= \mathbb C,\,\, \mathfrak v(\mathcal S)\simeq \mathfrak {pgl}_{k_1|l_1}(\mathbb C),\,\,
	(k_1,l_1)\ne (m,0)\,\, \text{and}\,\, (k_1,l_1)\ne (0,n),
	$$
	 then
	$\Ker\mathcal{P}= \{0\}$.}

\medskip

\noindent{\it Proof.}  
Consider the super-stabilizer $\mathcal H\subset \GL_{m|n}(\mathbb C)$ of $o$. It contains all super-matrices of the following form:
\begin{equation}\label{eq H contais matrices super case}
\left(
\begin{array}{cccc} 
A_1&0&* &0\\
C_1&B_1&*&D_1\\
*&0&A_2&0\\
*&D_2&C_2&B_2
\end{array} 
\right),
\end{equation} 
where the size of all matrices is as in Formula (\ref{eq H contais matrices}). Consider also the following Lie subsupergroup $\mathcal L$ in $\mathcal H$:
$$
\left(
\begin{array}{cc}
B_1 & D_1 \\
D_2 & B_2
\end{array}
\right).
$$
Clearly, $\mathcal L\simeq \GL_{k_1|l_1}(\mathbb C)$. Repeating computations (\ref{eq action over o}) for super-matrix (\ref{eq H contais matrices super case}), we see that $\mathcal L$ acts on $\mathcal S$ in the natural way, see (\ref{eq action of Q}), and the $\mathfrak l$-module $(\mathbf W_0)_o \simeq \mathfrak {pgl}_{k_1|l_1}(\mathbb C)$  is isomorphic to the adjoint $\mathfrak l$-module. Here $\mathfrak l\simeq \mathfrak {gl}_{k_1|l_1}(\mathbb C)$ is the Lie superalgebra of $\mathcal L$.

Let $\pi : \mathcal W\to \widetilde{\mathcal W}_0 = \mathcal W/\mathcal W_{(1)}$ be the natural map and $\pi_{o} : \mathcal W\to (\mathbf W_{0})_o$ be the composition of $\pi$ and of the evaluation map at the point $o$. We have the following commutative diagram: 
$$
\begin{CD}
\mathcal W(\mathcal B_0)@>{[X,\,\,\cdot
	\,\,]}>> \mathcal W(\mathcal B_0)
\\
@V{\pi_o}VV @V{\pi_o}VV\\
(\mathbf W_{0})_o@>{[X,\,\,\cdot
	\,\,]}>> (\mathbf W_{0})_o
\end{CD},
$$
where $X\in \mathfrak{l}$. (Note that the vector space $\mathcal W(\mathcal B_0)$ is an ideal in $\mathfrak v(\mathcal M)$ and in particular it is invariant with respect to the action of $\mathcal L$.)
Denote by $V$ the image $\pi_o(\mathcal W(\mathcal B_0))$. From the commutativity of this diagram it follows that 
$$
V\subset (\mathbf W_0)_o \simeq \mathfrak
{pgl}_{k_1|l_1}(\mathbb C)
$$
 is invariant with respect to the adjoint representation of $\mathfrak {pgl}_{k_1|l_1}(\mathbb C)$. Therefore, $V$ is an ideal in  $\mathfrak {pgl}_{k_1|l_1}(\mathbb C)$.

Let us describe ideals of the Lie superalgebra  $\mathfrak {pgl}_{k_1|l_1}(\mathbb C)$, where $(k_1,l_1)\ne (1,1)$, see \cite{Kac} for details. (The Lie superalgebra $\mathfrak {pgl}_{1|1}(\mathbb C)$ is nilpotent. We do not consider this case here because $\mathcal O_{\mathcal S}(\mathcal S_0)\ne \mathbb C$ for $\mathcal S= \mathbf {F}_{k'|l'}^{1|1}$.) This Lie superalgebra contains two trivial ideals $I=\{0\}$,
$\mathfrak {pgl}_{k_1|l_1}(\mathbb C)$ and it has one proper ideal 
$$
\mathfrak{psl}_{k_1|k_1}(\mathbb
C)=\mathfrak{sl}_{k_1|k_1}(\mathbb
C)/ \langle E_{2k_1}\rangle
$$
for $k_1=l_1$.

Clearly, we have $V\subset \Im(\gamma)$, where $\gamma : \widetilde{\mathcal W}_0(\mathcal B_0) \to (\mathbf W_{0})_o$ is the evaluation map. By Lemma \ref{lem sections of W_0}, we see that $\Im(\gamma)$ never coinsides with $\mathfrak {pgl}_{k_1|l_1}(\mathbb C)$ or $\mathfrak{psl}_{k_1|k_1}(\mathbb C)$. Hence, $V=\{0\}$.
In other words, all sections of  $\pi(\mathcal W(\mathcal B_0))$ are equal to $0$ at the point $o$. Since $\mathbf W_{0}$ is a homogeneous bundle, we get that $\pi(\mathcal W(\mathcal B_0))$ are equal to $0$ at any point. Therefore, $\pi(\mathcal W(\mathcal B_0))=\{0\}$ and 
$$
\mathcal W(\mathcal B_0)_{(0)} \simeq \mathcal W(\mathcal B_0)_{(1)}.
$$
 From Lemma  \ref{lem assume ker Pi ne 0} it follows that $\Ker\mathcal{P}=\{0\}$.$\Box$

\medskip

Using Theorem \ref{teor Ker Pi = 0} and Formula (\ref{eq connection mu and mu_B}), we get the following statement:

\medskip

\t\label{teor main generic case} {\sl Assume that $r>1$. If 
$$
\mathcal O_{\mathcal S}(\mathcal S_0)= \mathbb C,\,\, \mathfrak v(\mathbf F^{m|n}_{k_1|l_1})\simeq \mathfrak {pgl}_{m|n}(\mathbb C)\,\,\,\text{and} \,\,\, \mathfrak v(\mathbf
	F^{k_1|l_1}_{k'|l'})\simeq \mathfrak {pgl}_{k_1|l_1}(\mathbb C),
	$$ 
	then
	$$
	\mathfrak
	v(\mathbf F^{m|n}_{k|l})\simeq \mathfrak {pgl}_{m|n}(\mathbb C).
	$$
	}
\medskip

\subsection{Vector fields on flag supermanifolds, some exceptional cases}

\subsubsection{The base $\mathcal B$ is an exceptional super-Grassmannian}

 Assume that $r>1$,  $\mathcal O_{\mathcal S}(\mathcal S_0) = \mathbb C$ and $\mathcal B=\mathbf F_{k_1|l_1}^{m|n}$ 
 is one of the following super-Grassmanians:
 \begin{itemize}
 	\item [\textbf{a)}] $\mathbf F_{k_1|l_1}^{m|n} = \mathbf
 	F_{0|n}^{m|n}$ or $\mathbf F_{m|0}^{m|n}$, case 1.3 of  Theorem \ref{teor super=Grassmannians}.
 	
 	\item [\textbf{b)}] $\mathbf F_{k_1|l_1}^{m|n} = \mathbf F_{1|2}^{2|2}$ or $\mathbf F_{2|1}^{2|2}$, case 1.4 of Theorem \ref{teor super=Grassmannians}.  (We do not consider super-Grassmannians $\mathbf F_{1|0}^{2|2}$ and $\mathbf F_{0|1}^{2|2}$ here, because in these cases $\mathcal O_{\mathcal S}(\mathcal S_0)\ne \mathbb C$.)

 \item [\textbf{c)}] $\mathbf F_{k_1|l_1}^{m|n} =  \mathbf {F}^{1|n}_{0|n-1}$, where $n>2$, or $\mathbf {F}^{m|1}_{m-1|0}$, where $m>2$, case 1.2 of Theorem \ref{teor super=Grassmannians}. 
 This case we will consider in a separate paper.

  \item [\textbf{d)}] $\mathbf F_{k_1|l_1}^{m|n} =  \mathbf {F}^{2|2}_{1|1}$, case 1.1 of Theorem \ref{teor super=Grassmannians}. In this case $\mathcal O_{\mathcal S}(\mathcal S_0) \ne \mathbb C$. We do not consider this case here.
  
 \end{itemize}

\noindent{\bf Case \textbf{a}.} Without loss of generality we may consider only the case $\mathbf F_{k_1|l_1}^{m|n} = \mathbf
F_{0|n}^{m|n}$. In this case the base space $\mathbf F_{0|n}^{m|n}$ is a superpoint, i.e. it is a superdomain with the underlying space $\{\verb"pt"\}$,  one point, and with $mn$ odd coordinates.  Since $\mathbf
F_{k|l}^{m|n}$ is a superbundle with the base space isomorphic to a superpoint, we have 
$$
\mathbf
 F_{k|l}^{m|n}=\mathbf F_{0|n}^{m|n}\times \mathbf
 F_{k'|l'}^{0|n},\,\,\, \text{where}\,\,\, k'=(0,\ldots,0)\,\,\, \text{and}\,\,\,  l'=(l_2,\ldots,l_r).
 $$
 Our goal now is to prove the following theorem.
 \medskip
 
 \t \label{teor_W_nm+pgl_0n} {\it Assume that $r>1$ and $(k_1,l_1)=(m,0)$ or  $(k_1,l_1)=(0,n)$. Then 
 	$$
 	\mathfrak v(\mathbf F_{k|l}^{m|n})=W_{mn}\subset\!\!\!\!\!\!+
 	\Big(\bigwedge(mn)\otimes\mathfrak {pgl}_{n}(\mathbb C)\Big),
 	$$
 	where $W_{mn}=\operatorname{Der}(\bigwedge(mn))$. }
 
 \medskip
 
 \noindent {\it Proof.} The result follows from the following facts: 
\begin{align*}
\mathbf
F_{k|l}^{m|n}=\mathbf F_{0|n}^{m|n}\times \mathbf
F_{k'|l'}^{0|n},\quad  \mathcal O_{\mathcal S}(\mathcal S_0)=\mathbb C,\quad \mathcal O_{\mathcal B}(\mathcal B_0)= \bigwedge(mn),\\
 \mathfrak v(\mathbf F_{0|n}^{m|n}) \simeq W_{mn}, \quad \mathfrak v(\mathbf F_{0|l'}^{0|n}) \simeq  \mathfrak {pgl}_{n}(\mathbb C).
\end{align*}
In more details, since $ \mathcal O_{\mathcal S}(\mathcal S_0)=\mathbb C$, we have a Lie superalgebra homomorphism
$$
\mathcal P: \mathfrak v(\mathbf
F_{k|l}^{m|n})\to \mathfrak v(\mathbf F_{0|n}^{m|n}) \simeq W_{mn}.
$$
Since the bundle projection $\mathbf
F_{k|l}^{m|n} \to \mathbf F_{0|n}^{m|n}$ is just the projection to the first factor 
$$
\mathbf
F_{k|l}^{m|n}= \mathbf F_{0|n}^{m|n}\times \mathbf
F_{k'|l'}^{0|n}\to \mathbf F_{0|n}^{m|n},
$$
all vector fields on $\mathbf F_{0|n}^{m|n}$ can be lifted to $\mathbf
F_{k|l}^{m|n}$. The kernel of $\mathcal P$ is isomorphic to  $\bigwedge(mn)\otimes\mathfrak {pgl}_{n}(\mathbb C)$.
The proof is complete.$\Box$

 \bigskip

\noindent{\bf Case \textbf{b}.} Assume that $r=2$. Without loss of generality we may consider only the case $\mathbf F_{k_1|l_1}^{m|n} = \mathbf F^{2|2}_{1|2}$.  Under restriction $\mathcal O_{\mathcal S}(\mathcal S_0) = \mathbb C$ the fiber $\mathcal S$ can be one of the following super-Grassmanians: 
$$
\mathcal S =\mathbf F_{1|1}^{1|2}\,\,\,\text{or} \,\,\, \mathbf
F_{0|1}^{1|2}.
$$
 We have seen that $\mathfrak{v}(\mathbf
 F^{2|2}_{1|2})\simeq \tilde{\mathbf H}_4+\!\!\!\!\!\!\!\supset\langle
 z\rangle$, see (\ref{eq H_4 =}), Theorem \ref{teor super=Grassmannians}. A standard computation shows that  the image of $\mathfrak {gl}_{2|2}(\mathbb C)$ in $\mathfrak{v}(\mathbf
 F^{2|2}_{1|2})$ is
 $$
 (\tilde{\mathbf H}_4)_{-1}\oplus  (\tilde{\mathbf H}_4)_{0}\oplus  (\tilde{\mathbf H}_4)_{1}\oplus \langle z\rangle \simeq \mathfrak {pgl}_{2|2}(\mathbb C).
 $$ 
 Therefore,
\begin{equation}\label{eq excep case decomposition}
 \mathfrak{v}(\mathbf
 F^{2|2}_{1|2})\simeq \mathfrak {pgl}_{2|2}(\mathbb C)\oplus \langle \theta\rangle,
 \end{equation}
 as vector superspaces. (See (\ref{eq H_4 =}) for the definition of $\theta$.)  
 By Theorem \ref{teor constant functions} we have $\mathcal O_{\mathcal S}(\mathcal S_0)= \mathbb C$. Hence by Proposition \ref{prop Bash} we have a homomorphism of Lie superalgebras
 $$
 \mathcal{P}:\mathfrak{v}(\mathbf F^{2|2}_{k|l})\to \mathfrak
 v(\mathbf F^{2|2}_{1|2}).
 $$
  By Theorem \ref{teor super=Grassmannians} we see that $\mathfrak
 v(\mathcal S)\simeq \mathfrak {pgl}_{1|2}(\mathbb C)$. Therefore by Theorem \ref{teor Ker Pi = 0} the homomorphism $\mathcal{P}$
 is injective. The vector fields from $\mathfrak {pgl}_{2|2}(\mathbb C)$ are fundamental with respect to the action of the Lie superalgebra $\mathfrak {gl}_{2|2}(\mathbb C)$. Hence they can be lifted to the flag supermanifold $\mathbf F^{2|2}_{k|l}$. Therefore we need to find  $\mathcal{P}^{-1}(\theta)$. We will show that $\theta \notin \Im(\mathcal{P})$, i.e. $\theta$ cannot be lifted to $\mathbf F^{2|2}_{k|l}$.

 \medskip
 
 \t \label{t iskl_baza} {\it We have 
 	$$
 	\mathfrak{v}(
 	\mathbf F^{2|2}_{(1,1)|(2,1)})\simeq \mathfrak{pgl}_{2|2}(\mathbb C)\,\, \text{and}\,\,
 	\mathfrak{v}(\mathbf F^{2|2}_{(1,0)|(2,1)})\simeq
 	\mathfrak{pgl}_{2|2}(\mathbb C).
 	$$
 	}
 
 \medskip

\noindent{\it Proof.} Consider the following chart on $\mathbf F^{2|2}_{(1,1)|(2,1)}$:
 \begin{equation}
 	\label{fib_Z=} 
 	Z_{I_1}=\left(
 	\begin{array}{ccc}
 		x & \xi_1 & \xi_2 \\
 		1 & 0 & 0 \\
 		0 & 1 & 0 \\
 		0 & 0 & 1 \\
 	\end{array}
 	\right),\,\,Z_{I_2}= \left(
 	\begin{array}{cc}
 		1 & 0 \\
 		\eta & y \\
 		0 & 1 \\
 	\end{array}
 	\right)
 \end{equation}
Assume that $w:=\mathcal{P}^{-1}(\theta)$ is well-defined. Since all vector fields on $\mathbf F^{2|2}_{(1,1)|(2,1)}$ are projectable, in cootdinates (\ref{fib_Z=}) $w$ is equal to $\theta+v$, where $ v=f\frac{\partial}{\partial
 	y}+g\frac{\partial}{\partial \eta} $ is a vertical vector field and $f,\,g$ are holomorphic functions in coordinates (\ref{fib_Z=}). Let us find $f$ and $g$.
  We need the following fundamental vector fields on $\mathbf F^{2|2}_{(1,1)|(2,1)}$ written in coordinates (\ref{fib_Z=}):
 \begin{equation}\label{fund pol}
 \begin{split}
  E_{13} \longmapsto	\frac{\partial}{\partial \xi_1}, \quad  E_{14} \longmapsto
 	\frac{\partial}{\partial \xi_2},\quad  E_{42} \longmapsto \xi_2\frac{\partial}{\partial
 		x}+y\frac{\partial}{\partial \eta},\\
 	 E_{32} \longmapsto	\xi_1\frac{\partial}{\partial x}-\frac{\partial}{\partial \eta},
 	\quad  E_{34} \longmapsto -\xi_1\frac{\partial}{\partial \xi_2}-\frac{\partial}{\partial
 		y}.
 	\end{split}
 \end{equation}
 Here we denote by $E_{ij}$ the elementary matrix from $\mathfrak{gl}_{2|2}(\mathbb C)$.
 
  Since $\Ker\mathcal{P}=\{0\}$, using (\ref{fund pol}), we get
 \begin{align*}
& [\frac{\partial}{\partial \xi_1},
 w]=\xi_2\frac{\partial}{\partial x}+ \frac{\partial f}{\partial
 	\xi_1} \frac{\partial}{\partial y}+ \frac{\partial g}{\partial
 	\xi_1} \frac{\partial}{\partial \eta} = \xi_2\frac{\partial}{\partial
 	x}+y\frac{\partial}{\partial \eta};\\
&[\frac{\partial}{\partial \xi_2}, w]= -\xi_1\frac{\partial}{\partial
 	x}+ \frac{\partial f}{\partial \xi_2} \frac{\partial}{\partial y}+
 \frac{\partial g}{\partial \xi_2} \frac{\partial}{\partial \eta} =
 -\xi_1\frac{\partial}{\partial x}+\frac{\partial}{\partial \eta}.
  \end{align*}
 Hence, 
 $$
 \frac{\partial f}{\partial \xi_1}=0,\quad
 \frac{\partial g}{\partial \xi_1}=y,\quad
\frac{\partial f}{\partial \xi_2}=0, \quad 
 \frac{\partial g}{\partial \xi_2}=1.
 $$
 Furthermore,
 $$
 \begin{aligned}
 &[\xi_1\frac{\partial}{\partial \xi_2}+\frac{\partial}{\partial
 	y}, w]= \xi_1\frac{\partial}{\partial \eta}+ \frac{\partial
 	f}{\partial y} \frac{\partial}{\partial y}+ \frac{\partial
 	g}{\partial y}
 \frac{\partial}{\partial \eta} =0.
 \end{aligned}
 $$
 Hence, $\frac{\partial f}{\partial y}=0$ and $\frac{\partial g}{\partial y}= - \xi_1$. Now we see that
 $$
 \frac{\partial^2 g}{\partial \xi_1\partial y} = -1,\quad \frac{\partial^2 g}{\partial y \partial \xi_1} = 1.
 $$
 This is a contradiction. Therefore, $$
 \mathcal{P}^{-1}(z)=\emptyset\,\,\,
 \text{and}\,\,\, \mathfrak{v}(\mathbf F^{2|2}_{(1,1)|(2,1)})\simeq
 \mathfrak{pgl}_{2|2}(\mathbb C).
 $$
 
 The proof in the case $\mathbf
 F^{2|2}_{(1,0)|(2,1)}$ is similar.$\Box$
 
 \medskip

\subsubsection{The fiber $\mathcal S$ is an exceptional super-Grassmannian}

Assume that $r=2$,  $\mathcal O_{\mathcal S}(\mathcal S_0)=\mathbb C$ and $\mathcal S = \mathbf{F}^{k_1|l_1}_{k_2|l_2}$ is one of the following super-Grassmanians:
  \begin{itemize}
 	\item[\textbf{a)}] $\mathcal S=\mathbf{F}^{2|2}_{1|1}$, case $1.1$ of Theorem  \ref{teor super=Grassmannians};
 	\item [\textbf{b)}] $\mathcal S =\mathbf{F}^{2|2}_{0|1}$,
 	$\mathbf{F}^{2|2}_{1|0}$, $\mathbf{F}^{2|2}_{1|2}$ or $\mathbf{F}^{2|2}_{2|1}$, case $1.4$ of Theorem  \ref{teor super=Grassmannians};
 	\item[\textbf{c)}] $\mathcal S=\mathbf{F}^{1|l_1}_{0|l_1-1}$, $\mathbf{F}^{k_1|1}_{k_1-1|0}$, $\mathbf{F}^{1|l_1}_{1|1}$ or $\mathbf{F}^{k_1|1}_{1|1}$, where $n>2$,  case $1.2$ of Theorem  \ref{teor super=Grassmannians}.
 	\item[\textbf{d)}]  $\mathcal S = \mathbf{F}^{k_1|l_1}_{0|l_1}$ or $\mathbf{F}^{k_1|l_1}_{k_1|0}$, case 1.3 of Theorem  \ref{teor super=Grassmannians}. In both cases $\mathcal O_{\mathcal S}(\mathcal S_0)\ne \mathbb C$. We do not consider this case here.  		
  \end{itemize}

 Our goal now is to prove the following theorem.

 \medskip
 \t \label{teor fib_1} {\it Assume that $r=2$ and the fiber $\mathcal S$ of the superbundle $\mathbf{F}^{m|n}_{k|l}$ is a super-Grassmanian of type \textbf{a} or \textbf{b}. Then we have
$$
\mathfrak v(\mathbf{F}^{m|n}_{k|l})\simeq \mathfrak {pgl}_{m|n}(\mathbb C).
$$
 }
 
 \medskip
 
First of all let us compute the representation $\psi$ of the stabilizer $H$ in these cases. Formula (\ref{eq Z_I_2_preobrazuetsja}) tells us that the action of $H$ in $\mathcal S$ coinsides with the restriction of this action on $\GL_{2|2}(\mathbb C)_{\bar 0}$. We need the following lemma:
 
 \medskip
 
 \lem \label{lem fib_vesa}  {\sl The representation  $\psi$ of $H$ in the fiber $(\mathbf W_0)_o$ is completely reducible and its highest weights are:
 	\begin{enumerate}
 		\item $\mu_{m-1}-\mu_{m}$, $\lambda_{n-1}-\lambda_{n}$,
 		$\mu_{m-1}-\lambda_{n}$, $\lambda_{n-1}-\mu_{m}$, $0$,
 		$\mu_{m-1}+\mu_{m}-\lambda_{n-1}-\lambda_{n}$,
 		$\lambda_{n-1}+\lambda_{n}-\mu_{m-1}-\mu_{m}$, in case \textbf{a}.
 		\item $\mu_{m-1}-\mu_{m}$, $\lambda_{n-1}-\lambda_{n}$,
 		$\mu_{m-1}-\lambda_{n}$, $\lambda_{n-1}-\mu_{m}$, $0$,
 		$\mu_{m-1}+\mu_{m}-\lambda_{n-1}-\lambda_{n}$, in case \textbf{b}, super-Grassmannians $\mathbf F^{2|2}_{0|1}$ and $ \mathbf F^{2|2}_{1|2}$.
 		
 		\item $\mu_{m-1}-\mu_{m}$, $\lambda_{n-1}-\lambda_{n}$,
 		$\mu_{m-1}-\lambda_{n}$, $\lambda_{n-1}-\mu_{m}$, $0$,
 		$-\mu_{m-1}-\mu_{m}+\lambda_{n-1}+\lambda_{n}$, in case \textbf{b}, super-Grassmannians $\mathbf F^{2|2}_{1|0} $ and $\mathbf F^{2|2}_{2|1}$.
 	\end{enumerate}
 }
 
 \medskip
 
 \noindent{\it Proof.} As in Section $3.2$, we see that the nilradical of $H$ and the subgroup $\GL_{m-2}(\mathbb C)\times \GL_{n-2}(\mathbb C)$ in $H$ act trivialy on $\mathcal S$. The subgroup $\GL_{2}(\mathbb C)\times \GL_{2}(\mathbb C)$  acts in the natural way. Consider Decomposition (\ref{eq Gr_2|2,1|1 decomposition of Gl_0 mod}). We computed already highest weights of $\mathfrak{gl}_{2|2}(\mathbb C)_{\bar 0}$-module $\mathfrak{pgl}_{2|2}(\mathbb C)$. They are 
\begin{equation}\label{eq highest weigts Gr_2|2,1|1}
\mu_{m-1}-\mu_{m},\,\, \lambda_{n-1}-\lambda_{n},\,\,
\mu_{m-1}-\lambda_{n},\,\, \lambda_{n-1}-\mu_{m},\,\, 0.
\end{equation}
 Using the explicite description of $\mathfrak v(\mathbf{F}^{2|2}_{1|1})$ given by (\ref{eq vector fields GR_2|2,1|1 even}), (\ref{eq vector fields GR_2|2,1|1 odd}) and (\ref{eq vector fields GR_2|2,1|1 additional}), we get:
 \begin{align*}
 [\mu_{m-1}\mu(E_{11})+\mu_{m}\mu(E_{22})+\lambda_{n-1}\mu(E_{33})+\lambda_{n}\mu(E_{44}), \xi\frac{\partial}{\partial \eta}]=\\ (\mu_{m-1}+\mu_{m}-\lambda_{n-1}-\lambda_{n})\xi\frac{\partial}{\partial \eta};\\
 [\mu_{m-1}\mu(E_{11})+\mu_{m}\mu(E_{22})+\lambda_{n-1}\mu(E_{33})+\lambda_{n}\mu(E_{44}), \eta\frac{\partial}{\partial \xi}]=\\ (-\mu_{m-1}-\mu_{m}+\lambda_{n-1}+\lambda_{n})\eta\frac{\partial}{\partial \xi}.
 \end{align*}
 Here $E_{ii}$, where $i=1\ldots 4$, are elementary matrices from $\mathfrak{gl}_{2|2}(\mathbb C)_{\bar 0}$.
 The result follows.

 Let us prove the second statement. Consider $\mathbf F^{2|2}_{1|2}$ and decomposition (\ref{eq excep case decomposition}) of $\mathfrak{v}(\mathbf F^{2|2}_{1|2})$. We see easily that the vector subspaces $\langle \theta \rangle$ and $\mathfrak
 {pgl}_{2|2}(\mathbb C)$ are invariant with respect to the action of the Lie algebra  $\mathfrak
 {pgl}_{2|2}(\mathbb C)_{\bar 0}$. Again the vector space $\mathfrak {pgl}_{2|2}(\mathbb C)$ was decomposed into a sum of irreducible representations, see (\ref{eq submodules of pgl}). The highest weights of $\psi|\mathfrak {pgl}_{2|2}(\mathbb C)$ are given by (\ref{eq highest weigts Gr_2|2,1|1}). 
 Let us compute the highest weight of $\langle \theta \rangle$. 
 The image of the Cartan subalgebra 
 $$
 \operatorname{diag}(\mu_{m-1},\mu_{m})\times
 \operatorname{diag}(\lambda_{n - 1},\lambda_{n})
 $$
  with respect to the homomorphism 
 $\mu: \mathfrak {gl}_{2|2}(\mathbb C)_{\bar 0}
 \longrightarrow \mathfrak v(\mathbf F^{2|2}_{1|2})$ in chart (\ref{eq chart on F_2|2,1|2}) is given by
  \begin{align*}
 &\mu(E_{11})= x\frac{\partial}{\partial x}+\xi_1\frac{\partial}{\partial
 	\xi_1}+\xi_2\frac{\partial}{\partial \xi_2},\,\,\,
 \mu(E_{22}) = - x\frac{\partial}{\partial x},\\
 &\mu(E_{33}) =-\xi_1\frac{\partial}{\partial
 	\xi_1},\,\,\,
 \mu(E_{44}) =-\xi_2\frac{\partial}{\partial \xi_2}.
  \end{align*}
 We have
 \begin{align*} [\mu_{m-1}\mu(E_{11})+\mu_{m}\mu(E_{22})+\lambda_{n-1}\mu(E_{33})+\lambda_{n}\mu(E_{44}),
 \theta]=\\ (\mu_{m-1}+\mu_{m}-\lambda_{n-1}-\lambda_{n})\theta.
 \end{align*}
The result follows.

Computations in the cases $\mathbf F^{2|2}_{2|1}$, $\mathbf F^{2|2}_{0|1}$ and $\mathbf F^{2|2}_{1|0}$ are similar.$\Box$
  \medskip
 
 \noindent{\it Proof of Theorem \ref{teor fib_1}.} First of all let us compute the vector space of global sections of the vector bundle $\mathbf W_0$ using Theorem \ref{teor borel}.  The dominant highest weights of the representation $\psi$ are in case \textbf{a}:
 \begin{enumerate}
 	\item $0$ if $m>2$ and $n>2$;
 	\item $0$, $\mu_{1}-\mu_{2}$, $\mu_{1}-\lambda_{n}$,
 	$\mu_{1}+\mu_{2}-\lambda_{n-1}-\lambda_{n}$ for $m=2$ and $n>2$;
 	\item $0$, $\lambda_{1}-\lambda_{2}$, $\lambda_{1}-\mu_{m}$,
 	$\lambda_{1}+\lambda_{2}-\mu_{m-1}-\mu_{m}$ for $m>2$ and $n=2$.
 \end{enumerate}
 In case \textbf{b} for $\mathcal{O}_{\mathcal S}\simeq
 \mathbf{F}^{2|2}_{1|2}$ or $\mathbf{F}^{2|2}_{0|1}$ the dominant highest weights of
 $\psi$ are:
 \begin{enumerate}
 	\item $0$ for $m>2$, $n>2$;
 	\item $0$, $\mu_{1}-\mu_{2}$, $\mu_{1}-\lambda_{n}$,
 	$\mu_{1}+\mu_{2}-\lambda_{n-1}-\lambda_{n}$ for $m=2$, $n>2$;
 	\item $0$, $\lambda_{1}-\lambda_{2}$, $\lambda_{1}-\mu_{m}$,
 	for $m>2$, $n=2$.
 \end{enumerate}
 In case \textbf{b} for $\mathcal{O}_{\mathcal S}\simeq
 \mathbf{F}^{2|2}_{2|1}$ or $\mathbf{F}^{2|2}_{1|0}$ the dominant highest weights of
 $\psi$ are:
 \begin{enumerate}
 	\item $0$ for $m>2$, $n>2$;
 	\item $0$, $\mu_{1}-\mu_{2}$, $\mu_{1}-\lambda_{n}$ for $m=2$, $n>2$;
 	\item $0$, $\lambda_{1}-\lambda_{2}$, $\lambda_{1}-\mu_{m}$, $-\mu_{m-1}-\mu_{m}+\lambda_{1}+\lambda_{2}$
 	for $m>2$, $n=2$.
 \end{enumerate}

We restrict all weights on the Cartan subalgebra of 
 $\mathfrak{sl}_{m}(\mathbb C)\oplus\mathfrak{sl}_{n}(\mathbb C)
 \subset\mathfrak{gl}_{m}(\mathbb C)\oplus\mathfrak{gl}_{n}(\mathbb C)$.
By Theorem \ref{teor borel}, in case \textbf{a} we have:
 $$
\widetilde{\mathcal{W}}_0(\mathcal B_0)=\left\{
 \begin{array}{ll}
 \mathbb C, & m>2,\,n>2; \\
 \mathbb C\oplus \mathfrak{r}_1 \oplus
 \mathfrak{r}_2\oplus \mathfrak{r}_3, & m=2,\,n>2; \\
 \mathbb C\oplus \mathfrak{r}_4 \oplus
 \mathfrak{r}_5\oplus \mathfrak{r}_6, & m>2,\,n=2.
 \end{array}
 \right.
 $$
Without loss of generality we consider only the case \textbf{b}, $\mathcal{O}_{\mathcal S}\simeq
\mathbf{F}^{2|2}_{1|2}$ or $\mathbf{F}^{2|2}_{0|1}$. We have
 $$
\widetilde{\mathcal{W}}_0(\mathcal B_0)=\left\{
 \begin{array}{ll}
 \mathbb C, & m>2,\,n>2; \\
 \mathbb C\oplus \mathfrak{r}_1 \oplus
 \mathfrak{r}_2\oplus \mathfrak{r}_3, & m=2,\,n>2; \\
 \mathbb C\oplus \mathfrak{r}_4 \oplus
 \mathfrak{r}_5, & m>2,\,n=2.
 \end{array}
 \right.
 $$
Here $\mathfrak{r}_1$, $\mathfrak{r}_2$, $\mathfrak{r}_3$,
 $\mathfrak{r}_4$, $\mathfrak{r}_5$, $\mathfrak{r}_6$ are irreducible 
 $\mathfrak{sl}_{m}(\mathbb C)\oplus\mathfrak{sl}_{n}(\mathbb C)$-modules with highest weights $\mu_{1}-\mu_{2}$, $\mu_{1}-\lambda_{n}$,
 $\mu_{1}+\mu_{2}-\lambda_{n-1}-\lambda_{n}$,
 $\lambda_{1}-\lambda_{2}$, $\lambda_{1}-\mu_{m}$ and
 $\lambda_{1}+\lambda_{2}-\mu_{m-1}-\mu_{m}$, respectively, and $\mathbb C$
is the irreducible  $\mathfrak{sl}_{m}(\mathbb
 C)\oplus\mathfrak{sl}_{n}(\mathbb C)$-module with weight $0$.
 
 We use notations of Theorem \ref{teor Ker Pi = 0}. We have seen that $V$ is invariant with respect to the action of Lie superalgebra $\mathfrak{pgl}_{2|2}(\mathbb C)$. Consider the case  \textbf{a}. In case  $\widetilde{\mathcal{W}}_0(\mathcal B_0)=\mathbb C$, we have $V=\mathbb C$
 or $\{0\}$. Since $\mathfrak{pgl}_{2|2}(\mathbb C)$ does not have any $1$-dimensional ideals, the trivial module $\mathbb C$ is not 
 $\mathfrak{pgl}_{2|2}(\mathbb C)$-invariant. Hence, $V=\{0\}$. 
 Consider the case $\widetilde{\mathcal{W}}_0(\mathcal B_0)\simeq \mathbb C\oplus  \mathfrak{r}_1 \oplus \mathfrak{r}_2\oplus \mathfrak{r}_3$. As in Proof of Theorem  \ref{teor Ker Pi = 0}, we see that any combination of $H$-modules $\gamma(\mathbb C)$, $\gamma(\mathfrak{r}_1)$, $\gamma(\mathfrak{r}_2)$ and $\gamma(\mathfrak{r}_3)$ is not invariant with respect to $\mathfrak{pgl}_{2|2}(\mathbb C)$, see explicit description 
(\ref{eq vector fields GR_2|2,1|1 even}), (\ref{eq vector fields GR_2|2,1|1 odd}) and (\ref{eq vector fields GR_2|2,1|1 additional}).
Hence again $V=\{0\}$. We finish the proof similarly to Theorem \ref{teor Ker Pi = 0}.

Other cases are similar.$\Box$

 \medskip

\subsection{Main result}

We put $k_0=m$, $l_0=n$.

\medskip

\t \label{flag} {\sl Assume that $r>1$ and that we have the following restrictions on the flag type: 
	\begin{align*}
	&(k_i,l_i)\ne (k_{i-1},0),\,(0,l_{i-1}),\,\, i\geq 2;\\
		&(k_{i-1},k_i|l_{i-1},l_i)\ne (1,0|l_{i-1},l_{i-1}-1),\, 	(1,1|l_{i-1},1),\,\, i\geq 1;\\
			&(k_{i-1},k_i|l_{i-1},l_i)\ne (k_{i-1},k_{i-1}-1|1,0),\, 	(k_{i-1},1|1,1),\,\, i\geq 1;\\
		&k|l\ne 	(0,\ldots,0| n,l_2,\ldots,l_r),\,\, k|l\ne 	( m,k_2,\ldots,k_r|0,\ldots,0).
	\end{align*}
	 Then
		$$
		\mathfrak{v}(\mathbf{F}^{m|n}_{k|l})\simeq \mathfrak {pgl}_{m|n}(\mathbb
		C).
		$$
		If $k|l= (0,\ldots,0| n,l_2,\ldots,l_r)$ or $k|l= 	( m,k_2,\ldots,k_r|0,\ldots,0)$, then
		$$
		\mathfrak{v}(\mathbf{F}^{m|n}_{k|l})\simeq W_{mn}\subset\!\!\!\!\!\!+ (\bigwedge(\xi_1,\ldots,\xi_{mn})\otimes\mathfrak {pgl}_{n}(\mathbb
		C)),
		$$
		where $W_{mn}=\operatorname{Der}\bigwedge(\xi_1,\ldots,\xi_{mn})$.

}

\medskip

Note that the flag supermanifolds $\mathbf{F}^{m|n}_{k|l}$ and $\mathbf{F}^{n|m}_{l|k}$ are isomorphic.

\noindent{\it Elizaveta Vishnyakova}

\noindent {Max Planck Institute for Mathematics, Bonn}

\noindent{\emph{E-mail address:}
	\verb"VishnyakovaE@googlemail.com"}

\end{document}